\documentclass[11pt,fullpage]{article}
\usepackage{amsmath,amssymb,amsfonts,amsthm}
\newcommand{\be}{\begin{equation}}
\newcommand{\ee}{\end{equation}}
\newcommand{\bea}{\begin{eqnarray}}
\newcommand{\eea}{\end{eqnarray}}
\newcommand{\bean}{\begin{eqnarray*}}
\newcommand{\eean}{\end{eqnarray*}}
\newcommand{\ben}{\begin{equation}{nonumber}}
\newcommand{\een}{\end{equation}{nonumber}}

\numberwithin{equation}{section}

\newtheorem{dfn}{Definition}[section]
\newtheorem{thm}[dfn]{Theorem}
\newtheorem{lema}[dfn]{Lemma}
\newtheorem{pro}[dfn]{Proposition}
\newtheorem{coro}[dfn]{Corollary}
\newtheorem{xmpl}[dfn]{Example}
\newtheorem{rmrk}[dfn]{Remark}
\newcommand{\bdfn}{\begin{dfn}}
\newcommand{\bthm}{\begin{thm}}
\newcommand{\blema}{\begin{lema}}
\newcommand{\bpro}{\begin{pro}}
\newcommand{\bcoro}{\begin{coro}}
\newcommand{\bxmpl}{\begin{xmpl}}
\newcommand{\brmrk}{\begin{rmrk}}
\newcommand{\edfn}{\end{dfn}}
\newcommand{\ethm}{\end{thm}}
\newcommand{\elema}{\end{lema}}
\newcommand{\epro}{\end{pro}}
\newcommand{\ecoro}{\end{coro}}
\newcommand{\exmpl}{\end{xmpl}}
\newcommand{\ermrk}{\end{rmrk}}

\newcommand{\B}{\mathcal B}
\newcommand{\A}{\mathcal A}
\newcommand{\E}{\mathcal E}
\newcommand{\C}{\mathcal C}

\newcommand{\K}{\mathcal K}
\newcommand{\M}{\mathcal M}
\renewcommand{\H}{\mathcal H}

\newcommand{\half}{\frac{1}{2}}
\newcommand{\raro}{\rightarrow}
\newcommand{\lgl}{\langle}
\newcommand{\rgl}{\rangle}

\begin{document}
\begin{center}
{\bf{\large Quantum random walks and their convergence }}\\
{\large  Lingaraj Sahu {\footnote {This work is adopted from the Ph D thesis \cite{L}. The author is  grateful to Prof K. B. Sinha and  Dr Debashish Goswami for many  fruitful discussions during  the preparation of this and  would like to
acknowledge the support of National Board of Higher Mathematics, DAE, India and
 to the
DST-DAAD programme.}}}\\
{Stat-Math Unit, Indian Statistical Institute,\\ 203,
B.T. Road, Kolkata 700 108, India.
\\ email :
lingaraj\_r@isical.ac.in, lingaraj@gmail.com}
\end{center}
\begin{abstract}
Using coordinate-free basic operators on toy Fock spaces \cite{AP}, quantum random walks are defined  following the ideas in  \cite{LP,AP}. Strong convergence of   quantum random walks associated with bounded structure maps is proved under suitable assumptions, extendings the result obtained in  \cite{KBS} in case of one dimensional noise. To handle infinite dimensional noise we have used the coordinate-free language of quantum stochastic calculus developed in  \cite{GS1}.
\end{abstract}

\section{Toy  Fock spaces and  quantum random walks}
\subsection{Toy  Fock spaces}
Let ${\mathcal K} =L^2(\mathbb{R}_+,{\mathbf{k}}_0)$  where ${\mathbf{k}}_0$ is a Hilbert space with an orthonormal basis $\{e_i:i\ge 1\}.$ 
We note that, for any $n\ge 0,$ the $n$-fold symmetric tensor product
of ${\mathcal K} =L^2(\mathbb{R}_+,{\mathbf{k}}_0)$ and their direct sum can 
canonically be embedded in the symmetric Fock space $\Gamma:=\Gamma({\mathcal K}).$ 
 For any partition
$S\equiv(0=t_0<t_1<t_2\cdots)$ of $\mathbb{R}_+,$ the Fock  space $\Gamma ({\mathcal K})$
can be viewed as the infinite tensor product $\bigotimes_{n\ge 1}\Gamma_n$ of symmetric Fock  spaces
$\{\Gamma_n=\Gamma({\mathcal K}_{(t_{n-1},t_n]})\}_{n\ge 1}$ with 
respect to the stabilizing sequence 
$\Omega=\{ \Omega_n: n\ge 1\},$ where
$\Omega_n=\Omega_{(t_{n-1},t_n]}$ is the vacuum vector in $\Gamma_n.$\\

For any  $ 0\le s  \le t $ and $i\ge 1$ we define a vector $\chi_{(s,t]}^i:=\frac{1_{(s,t]}\otimes e_i}{\sqrt{t-s}}\in \mathcal
K_{(s,t]}.$ It is clear that
$\{\chi_{(s,t]}^i\}_{i\ge1}$ is an orthonormal
family in ${\mathcal K}_{(s,t]}$ and hence in $\Gamma_{(s,t]}.$ Here
we note that the Hilbert subspace $\mathbf {k_{(s,t]}}$ of
$\Gamma_{(s,t]}$ spanned  by these orthonormal vectors is canonically isomorphic to
${\mathbf{k}}_0.$ Let us consider the subspace $\hat{\mathbf k}_{(s,t]}={\mathbb {C}}\ \Omega_{(s,t]}\bigoplus \mathbf
{k_{(s,t]}}$ of $\Gamma$ and denote the space 
$\hat{\mathbf k}_{(t_{n-1},t_n]}$ by $ \hat{\mathbf k}_n,$ which is isomorphic to 
$ \hat{\mathbf{k}}_0:={\mathbb {C}}\bigoplus {\mathbf{k}}_0 .$
 Now we are in a position to define the toy Fock spaces.
 \bdfn
 The toy Fock space
associated with the partition $S$ of ${\mathbb R}_+$  is defined to be the subspace
$\Gamma(S):=\bigotimes_{n\ge 1} \hat{\mathbf k}_n$ with respect to
the stabilizing vector $\Omega=\otimes_{n\ge1} \Omega_n.$ \edfn
\noindent For notational simplicity we write $\chi_n^i$ for the vector
$\chi_{(t_{n-1}, t_n]}^i.$ Let
$\sqcap$ be the set of all finite subsets of $\mathbb{N} \times
\mathbb{N}.$ Thus an element $A \in
\sqcap$ is given by $A=\{m_1,i_1;m_2,i_2;\cdots m_n,i_n\}$ for some
$n$ with $ 1\le m_1<m_2 \cdots m_n < \infty.$ For  $A\in \sqcap,$ we
associate a vector
\[\chi_A=\Omega_1\otimes\Omega_2\otimes\cdots \chi_{m_1}^{i_1}\otimes\cdots
\chi_{m_2}^{i_2}\otimes\cdots
\chi_{m_n}^{i_n}\otimes\Omega_{m_n+1}\cdots
\] in the toy Fock space $\Gamma(S).$ Clearly this family
$\{\chi_A:A\in \sqcap\}$ forms an orthonormal basis for $\Gamma(S).$
Let $P(S)$ be the orthogonal projection of $\Gamma$onto the toy Fock
space $\Gamma(S).$ Without loss of generality  now onwards let us consider
toy Fock spaces $\Gamma(S_h)$ associated with  regular
partition  $S_h\equiv (0,h,\cdots)$ for some $h>0$ and denote the orthogonal projection by $P_h.$ The projection $P_h$ is given by\\
\[P_h=P_0 \oplus \mbox{$\bigoplus$}_{n\ge 1}\sum_{1\le m_1<m_2\cdots
<m_n}~\sum_{i_1,i_2\cdots i_n\ge 1} \mbox{$\bigotimes$}_{l=1}^n|\chi_{m_l}^{i_l}><\chi_{m_l}^{i_l}|,\] where
$P_0$ is the orthogonal projection of the symmetric Fock space $\Gamma$ onto the one dimensional Hilbert space
$\mathbb C \Omega.$
A simple computation shows that, for $f\in {\mathcal K},$ given by $f=\sum_{i\ge 1} f_i \otimes e_i$ with $f_i\in L^2(\mathbb R_+),$ 
\bean
&& P_h(\Omega)=\Omega,\\
&&P_hf=\sum_{m,i\ge 1} \frac{1}{\sqrt{h}}
\int_{(m-1)h}^{mh} f_i(s) ds ~\chi_m^i,\\
&& P_h \textbf{e}(f)=\Omega  \oplus \mbox{$\bigoplus$}_{n\ge 1} \frac{1}{\sqrt{n!}} \left [
\sum_{1\le m_1<m_2\cdots
<m_n}~\sum_{i_1,i_2\cdots i_n\ge 1} \mbox{$\bigotimes$}_{l=1}^n
\frac{1}{\sqrt{h}} \int_{(m_l-1)h}^{{m_l}h} f_{i_l}(s) ds~
\chi_{m_l}^{i_l}\right]
\eean
and furthermore,
\bean 
&& P_h\textbf{e}(f)=P_h\textbf{e}(f_{(k-1)h]})
P_h\textbf{e}(f_{[k]}) P_h\textbf{e}(f_{\left[kh\right .})~\mbox{and}\\
&& P_h\textbf{e}(f_{[k]})=\Omega_k \oplus \sum_{i\ge 1}
 \frac{1}{\sqrt{h}}
\int_{(k-1)h}^{kh} f_i(s) ds \frac{1_{((k-1)h,kh]} \otimes
e_i}{\sqrt{h}}.
\eean
Let us consider the  subspace $ {\mathcal {M}}$ of $L^2(\mathbb{R}_+,{\mathbf{k}}_0),$
given by
 \[{\mathcal {M}}=\{f\in L^2(\mathbb{R}_+,{\mathbf{k}}_0): f_i\in {\mathcal C}_c^1(\mathbb{R}_+)  
 ~\mbox{and}~ f_i=0 ~\mbox{for all but finitely many }~i\}.\]  Clearly
 ${\mathcal {M}}$ is a  dense subspace, so the algebraic tensor product 
 ${\mathbf{h}}_0\bigotimes \E(M)$ is dense in ${\mathbf{h}}_0\bigotimes \Gamma.$
For $f\in \M$ we define a constant $c_f:= \sum_{i\ge 1}\sup_{\tau} |f_i^\prime(\tau)|$ where  $f_i^\prime$ denotes the first derivative of the function $f_i.$
We have the following estimates which will be needed later. 
\begin{lema} \label{norm-diff}
For any $f\in \M, k\ge 1,$ \\
    $\|(1-P_h[k]) \textbf{e}(f_{[k]})\|\le h (c_f+\|f\|_\infty) \|\textbf{e}(f_{[k]})\|.$
\elema
\begin{proof}
(a). We have
\bean \lefteqn{\|(1-P_h[k]) \textbf{e}(f_{[k]})\|}\\
&=& \|(P_0+P_1-P_h)\textbf{e}(f_{[k]})+[1-P_0-P_1]\textbf{e}(f_{[k]})\|\\
&\le& \|f_{[k]}-P_hf_{[k]}\|+\|[1-P_0-P_1]\textbf{e}(f_{[k]})\|.
\eean
It is clear that $\|[1-P_0-P_1]\textbf{e}(f_{[k]})\|\le h\|f \|_\infty^2 \|\textbf{e}(f_{[k]})\|.$
Let us consider the first term,
\bean \lefteqn{\|f_{[k]}-P_hf_{[k]}\|^2}\\
&=& \sum_{i\ge 1} \|1_{[k]}(f_i-\frac{1}{h}\int_{[k]} f_i(s)ds)\|^2\\
&=& \sum_{i\ge 1} \int_{[k]}  dr|f_i(r)-\frac{1}{h}\int_{[k]} f_i(s)ds|^2 \\
&=& \sum_{i\ge 1} \frac{1}{h^2}\int_{[k]}  dr|\int_{[k]}(f_i(r)- f_i(s)ds|^2 \\
&\le & \sum_{i\ge 1} \frac{1}{h^2}\int_{[k]}  dr[\int_{[k]} h \sup |f_i^\prime(\tau)|ds]^2.
\eean
Since $c_f= \sum_{i\ge 1}\sup_{\tau} |f_i^\prime(\tau)|$, the  required estimate follows.
\end{proof}
By this lemma it can be proved that  the family of  othogonal projections $P_h$ converges strongly to identity operator
in $\Gamma$ as $h$ tends to $0.$ 
\subsection{Coordinate-free basic operators }
Here  we define  basic operators associated with toy Fock space
 $\Gamma(S_h)$ using the fundamental processes in coordinate-free language of quantum stochastic calculus, developed in \cite{GS1}  and obtained some useful estimates.  For $S\in \B({\mathbf{h}}_0), R\in \B({\mathbf{h}}_0,{\mathbf{h}}_0 \bigotimes {\mathbf{k}}_0)$ and
 $T\in \B({\mathbf{h}}_0 \bigotimes {\mathbf{k}}_0)$ let us define four basic operators as follows, for $k\ge 1,$
\begin{equation}
\begin{split}
&N_S^1[k]=S P_0[k]= P_0[k] \frac{\Lambda_S^1[k]}{h} ,\\
&N_R^2[k]= \frac{\Lambda_R^2[k]}{\sqrt{h}}P_1[k],\\
&N_R^3[k]=P_1[k] \frac{\Lambda_R^3[k]}{\sqrt{h}}, \\
&N_T^4[k]=P_1[k](\Lambda_T^4[k])P_1[k]P_h[k],
\end{split}
\end{equation}
where 
 \begin{equation}
 \begin {split}
&\Lambda_S^1[k]={\mathcal I}_S ((k-1)h,kh),\\
&\Lambda_R^2[k]= a_R((k-1)h,kh),\\
&\Lambda_R^3[k]=a^\dag_R((k-1)h,kh), \\
&\Lambda_T^4[k]=\Lambda_T((k-1)h,kh).\\
\end{split}
\end{equation}
For definition of coordinate-free fundamental processes ${\mathcal I}_S,a_R, a^\dag_R$ and $\Lambda_T$ we refer to \cite{GS1}.
All these maps  $ \B({\mathbf h}_0)\ni S\mapsto \Lambda_S^1[k],$
$ \B({\mathbf h}_0,{\mathbf h}_0 \bigotimes {\mathbf{k}}_0)\ni R\mapsto \Lambda_R^2[k],~\Lambda_R^3[k] $  and
$\B({\mathbf h}_0 \bigotimes {\mathbf{k}}_0)\ni T\mapsto \Lambda_T^4[k]$ are linear, and hence the maps    $ \B({\mathbf h}_0)\ni S\mapsto N_S^1[k], \B({\mathbf h}_0,{\mathbf h}_0 \bigotimes {\mathbf{k}}_0)\ni R\mapsto N_R^2[k],~N_R^3[k] $  and
$\B({\mathbf h}_0 \bigotimes {\mathbf{k}}_0)\ni T\mapsto N_T^4[k]$ are so.
 It is clear that the subspace $\Gamma(S_h)$ is invariant under  all these   operators $N^l$ and  their 
action  on ${\mathbf{h}}_0\otimes \Gamma:$ for $u\in {\mathbf{h}}_0, f\in L^2(\mathbb{R}_+,{\mathbf{k}}_0) $ are given by
\begin{equation} \label{Ne(f)}
\begin{split}
&N_S^1[k]u\textbf{e}(f_{[k]})=Su\otimes \Omega_{[k]} ,\\
&N_R^2[k]u\textbf{e}(f_{[k]})= \frac{\Lambda_R^2[k]}{\sqrt{h}}u \otimes f_{[k]}\\
&~~~~~~~~~~~~~~~~~~=\frac{1}{\sqrt{h}}\int_{[k]} R^* (uf(s))ds  \Omega_{[k]},\\
&N_R^3[k]u\textbf{e}(f_{[k]})= \frac{\Lambda_R^2[k]}{\sqrt{h}}
u \otimes\Omega_{[k]}  \\
&~~~~~~~~~~~~~~~~~~=\frac{1_{{\mathbf{h}}_0}\otimes 1_{[k]}}{\sqrt{h}} ~Ru \\
&N_T^4[k]u\textbf{e}(f_{[k]})=\Lambda_T^4[k])P_h[k]f_{[k]} \\
&~~~~~~~~~~~~~~~~~~= (1_{{\mathbf{h}}_0}\otimes 1_{[k]}) T u\otimes P_hf(\cdot).
\end{split}
\end{equation}
For any 
$S_1,S_2 \in  \B({\mathbf h}_0),
R_1,R_2 \in  \B({\mathbf h}_0,{\mathbf h}_0 \bigotimes {\mathbf{k}}_0)$
and $T_1,T_2 \in  \B({\mathbf h}_0 \bigotimes {\mathbf{k}}_0)$
we observe the following simple but useful identities, which are easy to derive 
\begin{itemize} \label{Nf-homo}
\item $(N_R^2[k])^2=(N_R^3[k])^2=0,$
~~~~~~ $N_{S_1}^1[k]~N_{S_2}^1[k]=N_{S_1 S_2}^1[k],$
\item $ N_{R_1}^2[k]~N_{R_2}^3[k]=N_{R_1^* R_2}^1[k],$
~~~~~~ $ N_S^1[k]~N_R^2[k]=N_{RS^*}^2[k],$
\item $ N_R^2[k]~N_T^4[k]=N_{T^*R}^2[k],$
~~~~~~~~~$ N_R^3[k]~N_S^1[k]=N_{RS}^3[k],$
\item $ N_T^4[k]~N_R^3[k]=N_{TR}^3[k],$
~~~~~~~~~ $ N_{R_1}^3[k]~N_{R_2}^2[k]=N_{R_1 R_2^*}^4[k],$
\item $ N_{T_1}^4[k]~N_{T_2}^4[k]=N_{T_1 T_2}^4[k],$
~~~~~~ $N_S^1[k]+N_{S\otimes 1_{{\mathbf{k}}_0}}^{4}[k]=S\otimes P_h[k].$
\end{itemize}
From  (\ref{Ne(f)} ) we have 
\begin{equation} \label{N-norm}
\begin{split}
&\|N_S^1[k]u\textbf{e}(f_{[k]})\|=\|Su\| ,\\ 
&\|N_R^2[k]u\textbf{e}(f_{[k]})\|\le \sqrt{h}\| R\|~\|u\|~\|f\|_\infty ,\\ 
&\|N_R^3[k]u\textbf{e}(f_{[k]})\|\le \|Ru\| \\ 
&\|N_T^4[k]u\textbf{e}(f_{[k]})\| \le \sqrt{h} \|T\| \| u\| \|f\|.
\end{split}
\end{equation}
Here we also note the following which can be verified easily  using Lemma 2.12 and Lemma 2.14 in  \cite{GS1},
\begin{equation} \label{Lambda-norm}
\begin{split}
&\|\Lambda_R^3[k]u\textbf{e}(f_{[k]})\|^2\\
&=  \|(1_{{\mathbf h}_0 \otimes 1_{[k]}})Ru\|^2 \|\textbf{e}(f_{[k]})\|^2 +\|\int_{[k]} R^*(uf(s)) ds\|^2 \|\textbf{e}(f_{[k]})\|^2,\\
&\|\Lambda_{T}^{4}[k]u\textbf{e}(f_{[k]})\|^2\\
&=\int_{[k]}\| T uf(s)\|^2  ds\|\textbf{e}(f_{[k]})\|^2
+\|\int_{[k]} \lgl f(s),T_{f(s)} \rgl~ ds~ u\textbf{e}(f_{[k]})\|^2.\\
\end{split}
\end{equation}
In the above expression $T_{f(s)}\in \B(\mathbf{h}_0,\mathbf{h}_0\bigotimes \mathbf{k}_0)$ and $\lgl f(s),T_{f(s)} \rgl\in \B(\mathbf{h}_0)$ define by
\[T_{f(s)}u=T u\otimes f(s),\forall u\in \mathbf{h}_0   ~\mbox{and}\]
\[\lgl \lgl f(s),T_{f(s)} \rgl  u,   v\rgl=\lgl T_{f(s)}u,v\otimes f(s) \rgl, \forall u,v\in \mathbf{h}_0.\]
For the  basic operators $N^l$'s we have the  estimates:
\blema \label{N-lambda}
(a). For any $k\ge 1$ and $u\in {\mathbf{h}}_0,f\in {\mathcal {M}},$ 
\begin{enumerate} 
\item $ \|\{h~N_S^1[k]-\Lambda_S^1[k]\}u\textbf{e}(f_{[k]})\|\le
 h^{\frac{3}{2}} \|Su\| \|f\|_{\infty} \|\textbf{e}(f_{[k]})\|,$\\
 \item $ \|\{\sqrt{h}~N_R^2[k]-\Lambda_R^2[k]\}u\textbf{e}(f_{[k]})\|\le h^{\frac{3}{2}} \|R\|~\|u\| \|f\|_{\infty}^2 \|\textbf{e}(f_{[k]})\|,$\\
 \item $ \|\{\sqrt{h}~N_R^3[k]-\Lambda_R^3[k]\}u\textbf{e}(f_{[k]})\|\le
 2 h \|Ru\| \|f\|_{\infty} \|\textbf{e}(f_{[k]})\|,$\\
 \item $ \|\{N_T^4[k]-\Lambda_T^4[k]\}u\textbf{e}(f_{[k]})\|\le
 2 h \|T\| (c_f+\|f\|_{\infty}^2) \|u\textbf{e}(f_{[k]})\|.$
   \end{enumerate}
(b). For any $k\ge 1$ and $u,v\in {\mathbf{h}}_0,f,g\in {\mathcal {M}},$ we have
\begin{enumerate} 
\item $ |\lgl v\textbf{e}(g_{[k]}), \{h~N_S^1[k]-\Lambda_S^1[k]\}u\textbf{e}(f_{[k]})\rgl|\\ \\
\le h^{\frac{3}{2}} \|Su\| \|f\|_{\infty} \|\textbf{e}(f_{[k]})\| \|v\textbf{e}(g_{[k]})\|,$ \\
 \item $ |\lgl v\textbf{e}(g_{[k]}), \{\sqrt{h}~N_R^2[k]-\Lambda_R^2[k]\}u\textbf{e}(f_{[k]})\rgl|\\ \\
 \le
 h^{\frac{3}{2}} \|R\|~\|u\| \|f\|_{\infty}^2 \|g\|_{\infty} \|\textbf{e}(f_{[k]})\| \|v\textbf{e}(g_{[k]})\|,$ \\ 
 \item $|\lgl v\textbf{e}(g_{[k]}), \{\sqrt{h}~N_R^3[k]-\Lambda_R^3[k]\}u\textbf{e}(f_{[k]})\rgl |\\ \\
 \le
 2 h^2 \|Ru\|~\|v\|~\|f\|_{\infty}~ \|g\|_{\infty}~\|\textbf{e}(f_{[k]})\|^2~ \|\textbf{e}(g_{[k]})\|^2,$ \\
 \item $ |\lgl v\textbf{e}(g_{[k]}), \{N_T^4[k]-\Lambda_T^4[k]\}u\textbf{e}(f_{[k]})\rgl |\\ \\
 \le
 h^2 \left [(\|f\|_{\infty}+c_f)\|g\|_{\infty}\right ]^2  \|T\|~\|u\|~\|v\|~ \|\textbf{e}(f_{[k]})\|^2 ~\|\textbf{e}(g_{[k]})\|^2.$
   \end{enumerate}
 \elema

\begin{proof}
\noindent a.(1) It is clear from the definition that  \bean &&\|\{h~N_S^1[k]-\Lambda_S^1[k]\}u\textbf{e}(f_{[k]})\|
= h\|Su(\Omega_{[k]}-\textbf{e}(f_{[k]}))\|\\
&&~~~~~~~~~~~~~~~~~~~~~~~~~~~~~~~~=h\|Su\| \|\Omega_{[k]}-\textbf{e}(f_{[k]})\|\\
&&~~~~~~~~~~~~~~~~~~~~~~~~~~~~~~~~\le h^{\frac{3}{2}} \|Su\| \|f\|_{\infty} \|\textbf{e}(f_{[k]})\|.
\eean
\noindent (2) From  the definitions, we have \bean && \|\{\sqrt{h}~N_R^2[k]-\Lambda_R^2[k]\}u\textbf{e}(f_{[k]})\|=\|\int_{[k]} R^*(uf(s)) ds~ (\Omega_{[k]}-\textbf{e}(f_{[k]}))\|\\
&&~~~~~~~~~~~~~~~~~~~~~~~~~~~~~~~~~~
\le \int_{[k]} \|R^*(uf(s))\| ds~ \|(\Omega_{[k]}-\textbf{e}(f_{[k]}))\|\\
&&~~~~~~~~~~~~~~~~~~~~~~~~~~~~~~~~~~\le h^{\frac{3}{2}} \|R\|~\|u\| \|f\|_{\infty}^2 \|\textbf{e}(f_{[k]})\|. 
\eean
\noindent(3) We have
\bean\lefteqn
{\|\{\sqrt{h}~N_R^3[k]-\Lambda_R^3[k]\}u\textbf{e}(f_{[k]})\|^2}\\
&=&\|(1_{{\mathbf{h}}_0}\otimes 1_{[k]})~Ru - \Lambda_R^3[k]u\textbf{e}(f_{[k]})\|^2\\
&=& \|(1_{{\mathbf{h}}_0}\otimes 1_{[k]}) ~Ru\|^2 + \|\Lambda_R^3[k]u\textbf{e}(f_{[k]})\|^2\\
&&~~- 2 \mathcal{R}e \lgl (1_{{\mathbf{h}}_0}\otimes 1_{[k]}) ~Ru , \Lambda_R^3[k]u\textbf{e}(f_{[k]}) \rgl
\eean
Now using (\ref{Lambda-norm}) and the definition of $\Lambda_R^3$  the above quantity is equal to
\bean
&&  \|(1_{{\mathbf{h}}_0}\otimes 1_{[k]})Ru\|^2 +  \|(1_{{\mathbf{h}}_0}\otimes 1_{[k]})Ru\|^2 \|\textbf{e}(f_{[k]})\|^2
\\
&&~~+\|\int_{[k]} R^*(uf(s)) ds\|^2 \|\textbf{e}(f_{[k]})\|^2- 2  \|(1_{{\mathbf{h}}_0}\otimes 1_{[k]})Ru\|^2\\
&&=   \|(1_{{\mathbf{h}}_0}\otimes 1_{[k]})Ru\|^2 [\|\textbf{e}(f_{[k]})\|^2-1]
+\|\int_{[k]} R^*(uf(s)) ds\|^2 \|\textbf{e}(f_{[k]})\|^2\\
&&\le 
 2 h^2  \|R\|^2 \|u\|^2 \|f\|_{\infty}^2 \|\textbf{e}(f_{[k]})\|^2.
 \eean
\noindent(4) We have
 \bean 
 &&\|\{N_T^4[k]-\Lambda_T^4[k]\}u\textbf{e}(f_{[k]})\|^2\\
 &&= \|(1_{{\mathbf{h}}_0}\otimes 1_{[k]}) T (u\otimes P_hf(\cdot))\|^2+
\|\Lambda_T^4[k]u\textbf{e}(f_{[k]})\|^2\\
&&~~- 2 \mathcal{R}e \lgl  (1_{{\mathbf{h}}_0}\otimes 1_{[k]}) T (u\otimes P_hf(\cdot)),
\Lambda_T^4[k]u\textbf{e}(f_{[k]}) \rgl.
\eean
By the  definition of $\Lambda_T^4$ (see \cite{GS1})
\bean
&&\lgl  (1_{{\mathbf{h}}_0}\otimes 1_{[k]}) T (u\otimes P_hf(\cdot)),
\Lambda_T^4[k]u\textbf{e}(f_{[k]}) \rgl\\
&&=\lgl (1_{{\mathbf{h}}_0}\otimes 1_{[k]}) T (u\otimes P_hf(\cdot)),
 a^\dag (T_{f_{[k]}}^{[k]})u\textbf{e}(f_{[k]}) \rgl\\
&&=\lgl (1_{{\mathbf{h}}_0}\otimes 1_{[k]}) T (u\otimes P_hf(\cdot)),
 T_{f_{[k]}}^{[k]}
(u\Omega_{[k]}) \rgl\\
&&=\int_{[k]}\lgl T (uP_h(f)(s)),T( uf(s))\rgl~ ds.
\eean
Thus using  (\ref{Lambda-norm}) we obtained
\bean
&&\|\{N_T^4[k]-\Lambda_T^4[k]\}u\textbf{e}(f_{[k]})\|^2\\
&&=\int_{[k]}\| T (uP_h(f)(s))\|^2  ds\\
&&~~+\int_{[k]}\| T ( uf(s))\|^2  ds\|\textbf{e}(f_{[k]})\|^2
+\|\int_{[k]} \lgl f(s),T_{f(s)} \rgl~ ds~ u\textbf{e}(f_{[k]})\|^2\\
&&~~-2 \mathcal{R}e \int_{[k]}\lgl T (uP_h(f)(s)),T( uf(s))\rgl ds\\
&&=\int_{[k]}\| T uf(s)\|^2  ds(\|\textbf{e}(f_{[k]})\|^2-1)
+\|\int_{[k]} \lgl f(s),T_{f(s)} \rgl~ ds~ u\textbf{e}(f_{[k]})\|^2\\
&&~~+\int_{[k]}\| T (u\otimes (1-P_h)(f)(s))\|^2  ds\\
&& \le 2 h^2 \|T\|^2  \|f\|_{\infty}^4 \|u\textbf{e}(f_{[k]})\|^2+ \| T\|^2 \|u\|^2 \int_{[k]} \|(1-P_h)(f)(s))\|^2  ds  \\
&& \le 2  h^2 \|T\|^2  \|f\|_{\infty}^4 \|u\textbf{e}(f_{[k]})\|^2+ \| T\|^2 \|u\|^2  \|(1-P_h)(f_{[k]})\|^2 .
\eean
Since $\|(1-P_h) \textbf{e}(f_{[k]})\|^2 \le  h^2 c_f,$ the required estimate follows.

\noindent (b). The estimates (1) and (2) follow directly from (a).

\noindent (3) From the definitions
 \bean\lefteqn
{\lgl v\textbf{e}(g_{[k]}), \{\sqrt{h}~N_R^3[k]-\Lambda_R^3[k]\}u\textbf{e}(f_{[k]})\rgl }\\
&=&\lgl v\textbf{e}(g_{[k]}), \sqrt{h}~N_R^3[k]u\textbf{e}(f_{[k]})\rgl -\lgl v\textbf{e}(g_{[k]}), \Lambda_R^3[k] u\textbf{e}(f_{[k]})\rgl\\
&=&\lgl v\textbf{e}(g_{[k]}), (1_{{\mathbf{h}}_0}\otimes 1_{[k]})~Ru \rgl-
\lgl \Lambda_R^2[k] v\textbf{e}(g_{[k]}), u\textbf{e}(f_{[k]})\rgl\\
&=& \int_{[k]} \lgl R u,  vg(s)\rgl ds (1-\lgl \textbf{e}(g_{[k]}), \textbf{e}(f_{[k]})\rgl ). 
\eean
Thus we have obtained the required estimate,
\bean\lefteqn
{|\lgl v\textbf{e}(g_{[k]}), \{\sqrt{h}~N_R^3[k]-\Lambda_R^3[k]\}u\textbf{e}(f_{[k]})\rgl| }\\ 
&\le &h^2 \|Ru\| ~ \|v\|~\|f\|_{\infty} \|g\|_{\infty}\|^2\textbf{e}(f_{[k]})\|^2 \|\textbf{e}(g_{[k]})\|^2.
 \eean
 4. By definition of $N_T^4$ and $\Lambda_T^4$
  \bean\lefteqn
{\lgl v\textbf{e}(g_{[k]}), \{N_T^4[k]-\Lambda_T^4[k]\}u\textbf{e}(f_{[k]})\rgl }\\
&=&\lgl v\textbf{e}(g_{[k]}), N_T^4[k]u\textbf{e}(f_{[k]})\rgl -\lgl v\textbf{e}(g_{[k]}), \Lambda_T^4[k] u\textbf{e}(f_{[k]})\rgl
\eean
\bean
&=&\lgl v\textbf{e}(g_{[k]}), (1_{{\mathbf{h}}_0}\otimes 1_{[k]}) T(uP_hf(\cdot)) \rgl-
\lgl  v\textbf{e}(g_{[k]}), a^\dag (T_{f_{[k]}}^{[k]})u\textbf{e}(f_{[k]})\rgl\\
&=& \int_{[k]} \lgl vg(s),T(u(P_hf)(s))\rgl ds - \int_{[k]} \lgl vg(s),T(uf(s))\rgl ds \lgl \textbf{e}(g_{[k]}), \textbf{e}(f_{[k]})\rgl  . \\
&=& \int_{[k]} \lgl vg(s),T[u\left ((P_h-1)f\right )(s)]\rgl ds \\
&&+ \int_{[k]} \lgl vg(s),T(uf(s)) \rgl ds[1- \lgl \textbf{e}(g_{[k]}), \textbf{e}(f_{[k]})\rgl] . 
\eean
So we get
\bean\lefteqn
{|\lgl v\textbf{e}(g_{[k]}), \{N_T^4[k]-\Lambda_T^4[k]\}u\textbf{e}(f_{[k]})\rgl |}\\
&\le & \left (\int_{[k]} \|vg(s)\|^2  ds \right )^{\frac{1}{2}}  \left (\int_{[k]} \|T[u((P_h-1)f_{[k]}(s))]\|^2  ds \right )^{\frac{1}{2}} \\
&&+ \int_{[k]} \|vg(s)\| \| T(uf(s))\|   ds~ \|1- \lgl \textbf{e}(g_{[k]}), \textbf{e}(f_{[k]})\rgl\| \\
&\le & h \|v\|\| g\|_\infty  \|T\| \|u\| \|(P_h[k]-1)f_{[k]}\| \\
&&+ h\|v \| \| g\|_\infty \| T\| \| u\| \| f\|_\infty~ \|1- \lgl \textbf{e}(g_{[k]}), \textbf{e}(f_{[k]})\rgl\| .
\eean
Using the estimates of $\|(P_h[k]-1)f_{[k]}\| $ and $\|1- \lgl \textbf{e}(g_{[k]}), \textbf{e}(f_{[k]})\rgl\|$ the required estimate follows.
\end{proof}
\brmrk
The  estimates in the above Lemma will also hold if we replace the initial Hilbert space ${\mathbf h}_0$ by ${\mathbf h}_0 \bigotimes \Gamma_{\left .(k-1)h\right ]}$ and  take  $S\in \B({\mathbf h}_0 \bigotimes \Gamma_{\left .(k-1)h\right ]}),$\\$ R\in \B({\mathbf h}_0 \bigotimes \Gamma_{\left .(k-1)h\right ]},{\mathbf h}_0 \bigotimes \Gamma_{\left .(k-1)h\right ]}  \bigotimes {\mathbf{k}}_0)$ and
 $T\in \B({\mathbf h}_0 \bigotimes \Gamma_{\left .(k-1)h\right ]} \bigotimes {\mathbf{k}}_0).$
\ermrk
\subsection{Quantum random walk}
Let $\A \subseteq \B({\mathbf{h}}_0)$ be a von Neumann algebra. Let us   consider the Hilbert von Neumann  module $\A \bigotimes {\mathbf{k}}_0.$ 
Suppose we are given with a family of $*$-homomorphisms $\{\beta(h)\}_{h > 0}$    from $\A$ to
 $\A\bigotimes\B(\hat{\mathbf{k}}_0).$ For $h>0, \beta(h)$
can be written as\\ $\beta(h,x)=\left( \begin{array}{cc} \beta_1(h,x) & (\beta_2(h,x))^*\\
\beta_3(h,x) &  \beta_4(h,x)
\end{array} \right), \forall x\in \A,$ where the components $\beta_l(h)$'s are contractive  maps  and $\beta_1(h)\in \B(\A),
 \beta_4(h)\in \B(\A, \A\bigotimes \B({\mathbf{k}}_0))$ 
 and $\beta_2(h), \beta_3(h) \in \B(\A,\A\bigotimes {\mathbf{k}}_0).$ 
The $*$-homomorphic properties of  $\beta(h)$ can be translated into the following properties of $\beta_l(h)$'s.
\begin{itemize} 
\item $\beta_1(h,x^*)=(\beta_1(h,x))^*,~~\beta_4(h,x^*)=(\beta_4(h,x))^*,~~\beta_3(h,x^*)=\beta_2(h,x),$
\item $\beta_1(h,xy)=\beta_1(h,x)\beta_1(h,y)+(\beta_2(h,x))^*\beta_3(h,y),$
\item $\beta_2(h,xy)=\beta_1(h,x)(\beta_2(h,y))^*+(\beta_2(h,x))^*\beta_4(h,y),$
\item $\beta_3(h,xy)=\beta_3(h,x)\beta_1(h,y)+\beta_4(h,x)\beta_3(h,y),$
\item $\beta_4(h,xy)=\beta_3(h,x)(\beta_2(h,y))^*+\beta_4(h,x)\beta_4(h,y).$
\end{itemize}
We define a family of  maps ${\mathcal P}_t^{(h)}:\A \bigotimes \E(\K)\raro \A\bigotimes \Gamma$ 
as follows. We subdivide the interval $[0,t]$ 
into $[k]\equiv \left( (k-1)h,kh\right],~1\le k\le n$ so that $t\in \left( (n-1)h,nh\right]$ as earlier  and set for $x\in \A,~f\in \K$
\begin{equation}
\left.
\begin{array}{l} 
 {\mathcal P}_{0}^{(h)}(x \textbf{e}(f))=x \textbf{e}(f)\\ \\
{\mathcal P}_{kh}^{(h)}(x \textbf{e}(f))=\sum_{l=1}^4
{\mathcal P}_{(k-1)h}^{(h)}  N_{\beta_l(h,x)}^l[k]  \textbf{e}(f)
\end{array}
\right \}
\end{equation}
and 
${\mathcal P}_{t}^{(h)}={\mathcal P}_{nh}^{(h)}.$\\
Now setting a  family of linear maps $p_t^{(h)}:\A\raro \A\bigotimes \B(\Gamma),$ by\\
 $p_t^{(h)}(x)u \textbf{e}(f):={\mathcal P}_t^{(h)}(x \textbf{e}(f))u, \forall u\in {\mathbf h}_0$ we have
\begin{equation}
\left.
\begin{array}{l}
 p_0^{(h)}(x)u\textbf{e}(f)=xu\textbf{e}(f)\\ \\
p_{t}^{(h)}(x)u \textbf{e}(f)= p_{nh}^{(h)}(x)u\textbf{e}(f)=\sum_{l=1}^4
 N_{p_{(n-1)h}^{(h)} (\beta_l(h,x))}^l[n]u\textbf{e}(f).
\end{array}
\right \}
\end{equation}
As per our convention  $p_{(n-1)h}^{(h)} $ appear above  are identified with their  ampliations 
$p_{(n-1)h}^{(h)} \otimes 1_{{\mathbf{k}}_0}$ as well as  $p_{(n-1)h}^{(h)} \otimes 1_{\B({\mathbf{k}}_0)}.$
For $k\ge 1, l=1,2,3$ and $4,  N_{p_{(k-1)h}^{(h)} (\beta_l(h,x))}^l[k]$ are  defined in terms of  
\[ \Lambda_{(p_{(k-1)h}^{(h)} (\beta_1(h,x))}^1[k],~
\Lambda_{(p_{(k-1)h}^{(h)} \otimes 1_{{\mathbf{k}}_0}) (\beta_2(h,x))}^2[k],
~\Lambda_{(p_{(k-1)h}^{(h)} \otimes 1_{{\mathbf{k}}_0}) (\beta_3(h,x))}^3[k]\]
and $\Lambda_{(p_{(k-1)h}^{(h)} \otimes 1_{\B({\mathbf{k}}_0)}) (\beta_4(h,x))}^4[k] $
where,  for example  
$\Lambda_{(p_{(k-1)h}^{(h)}\otimes 1_{{\mathbf{k}}_0}) (\beta_2(h,x))}^2[k]$ carries the  meaning of $a^\dag_{(p_{(k-1)h}^{(h)}\otimes 1_{{\mathbf{k}}_0}) (\beta_2(h,x))}[k]$ with initial Hilbert space ${\mathbf{h}}_0 \bigotimes \Gamma_{(k-1)h]}.$

\noindent For notational simplicity, for any bounded $*$-preserving  map \[\alpha: \A\raro \A\bigotimes\B(\hat{\mathbf{k}}_0),~
\alpha(x)=\left( \begin{array}{cc} \alpha_1(h,x) & ( \alpha_2(h,x))^*\\
 \alpha_3(h,x) &   \alpha_4(h,x)
\end{array} \right),\]   we write
$N_{ \alpha(h,x)}[k] $ for  $ \sum_{l=1}^4 N_{ \alpha_l(h,x)}^l[k]$ and
$\Lambda_{ \alpha(h,x)}[k] $ for $ \sum_{l=1}^4 \Lambda_{ \alpha_l(h,x)}^l[k].$
Now for each $k\ge 1$ defining a linear map  $\rho_k(h,x)=
N_{\beta(h,x)}[k], p_{nh}^{(h)}$  can be written as
$p_{nh}^{(h)}=\rho_1(h)\cdots \rho_n(h).$
By the properties of the family $\{\beta_l(h)\}$ and $\{N^l[k]\},$  each $\rho_k(h)$ is a $*$-homomorphism and hence $p_t^{(h)}$ is so. We call this family of homomorphisms $\{p_t^{(h)}:t\ge 0\}$ a {\it quantum random walk}.
In the next section we shall construct  quantum random walks associated with 
uniformly continuous QDS on a  von Neunmann algebra and  show the strong  convergence. Let us conclude this section 
with  the following observation which will be needed latter.

\blema \label{JPh} For any $t\ge 0, t\in ((n-1)h,nh]$ for some $n\ge 1$ and $ x\in \A,u\in {\mathbf h}_0$ and $f\in \K$\\
\be \label{jh-k} {\mathcal P}_t^{(h)}(x\textbf{e}(f))u= xu\textbf{e}(f)+\sum_{k=1}^n
{\mathcal P}_{(k-1)h}^{(h)}  N_{{\beta(h,x)}- b(x)}[k]\textbf{e}(f)u+ F(h,x,u,f),
\ee
 where $b(x)=\left( \begin{array}{cc} b_1(x) & (b_2(x))^*\\
b_3(x) & b_4(x) \end{array}\right)=x\otimes 1_{\hat{\mathbf{k}}_0}$  
and \\
 $F(h,x,u,f)=-\sum_{k=1}^n {\mathcal P}_{(k-1)h}^{(h)}  (x(1_{\Gamma}-P_h[k])\textbf{e}(f))u.$ Moreover,
 for any $f\in \M$
\be \label{F}\|F(h,x,u,f)\|^2\le h~ c(f,t) \|x\|^2~ \|u\|^2,\ee
where $c(f,t)=2t(c_f+\|f\|_\infty)\| \textbf{e}(f)\|.$
\elema
  \begin{proof}
Since for any $k\ge 1,$
\[N_{b(x)}[k]=\sum_{l=1}^4 N_{b_l(x)}^l[k]=N_x^1[k]+N_{x\otimes 1_{{\mathbf{k}}_0}}^4[k]=x\otimes P_h[k],\] 
We get
\bean
 \lefteqn
  { {\mathcal P}_t^{(h)}(x\textbf{e}(f))u={\mathcal P}_{nh}^{(h)}(x\textbf{e}(f))u}\\
&=& xu\textbf{e}(f)+\sum_{k=1}^n
({\mathcal P}_{kh}^{(h)}-{\mathcal P}_{(k-1)h}^{(h)})(x\textbf{e}(f))u\\
&=& xu\textbf{e}(f)+\sum_{k=1}^n
{\mathcal P}_{(k-1)h}^{(h)}  N_{{\beta(h,x)}- b(x)}[k]\textbf{e}(f)u
\eean
\bean
&&~~~~~~~~~~~~~~~-\sum_{k=1}^n
{\mathcal P}_{(k-1)h}^{(h)} (x\otimes 1_{\Gamma}- N_{ b(x)}[k])\textbf{e}(f)u\\
&&~~~~~~~~~~~~~= xu\textbf{e}(f)+\sum_{k=1}^n
{\mathcal P}_{(k-1)h}^{(h)}  N_{{\beta(h,x)}- b(x)}[k]\textbf{e}(f)u+F(h,x,u,f).
\eean
In order to obtained  ~(\ref{F}) let us consider the following. For any $1\le m\le n$ setting $ Z_{m}=\sum_{k=1}^m p_{(k-1)h}^{(h)}  (x)(1-P_h[k]),$ we have \\
$\| Z_{m}u\textbf{e}(f_{\left . mh\right ]})\| 
\le \sum_{k=1}^m \|p_{(k-1)h}^{(h)}  (x)u\textbf{e}(f_{(k-1)h]})\| ~ \| (1-P_h[k])\textbf{e}(f_{[k]})\| \|\textbf{e}(f_{(kh,mh]})\|.
$ Now using  Lemma ~\ref{norm-diff}(a)  and the fact that $p_{kh}^{(h)} $'s are homomorphisms,
 \bean &&\| Z_{m}u\textbf{e}(f_{\left . mh\right ]})\|\\
 &&\le \sum_{k=1}^m h (c_f+\|f\|_\infty) \|x\| \|u\textbf{e}(f_{\left . mh\right ]})\|\\
 &&\le t(c_f+\|f\|_\infty) \|x\| \|u\textbf{e}(f_{\left . mh\right ]})\|.
 \eean
We have
\bean && \|F(h,x,u,f)\|^2\\
&&=\sum_{k=1}^n \|p_{(k-1)h}^{(h)} (x) u\textbf{e}(f_{(k-1)h]})\|^2 \|(1-P_h[k])\textbf{e}(f_{[k]})\|^2 ~\|\textbf{e}(f_{\left[kh\right .})\|^2\\
&&~~+2 \mathcal{R}e \sum_{k=1}^n \lgl  Z_{k-1} u\textbf{e}(f_{(k-1)h]})  , p_{(k-1)h}^{(h)}  (x)   u\textbf{e}(f_{(k-1)h]})  \rgl\\
&& ~~~~~~~~~~~~~~~~~~~~~~\lgl \textbf{e}(f_{[k]}), (1-P_h[k])\textbf{e}(f_{[k]})\rgl  ~\|\textbf{e}(f_{\left[kh\right .})\|^2\\
&&\le \sum_{k=1}^n \|x\|^2 \|u\textbf{e}(f_{(k-1)h]})\|^2 \|(1-P_h[k])\textbf{e}(f_{[k]})\|^2 ~\|\textbf{e}(f_{\left[kh\right .})\|^2\\
&&~~+2 \sum_{k=1}^n   \|Z_{k-1} u\textbf{e}(f_{(k-1)h]})\|  ~\|x\| \|u\textbf{e}(f_{(k-1)h]})\| ~\\
&&~~~~~~~~~~~~~~ \|(1-P_h[k])\textbf{e}(f_{[k]})\|^2 ~\|\textbf{e}(f_{\left[kh\right .})\|^2.
\eean
Using the uniform bound for $\|Z_{k-1} u\textbf{e}(f_{(k-1)h]}) \|$ and  Lemma ~\ref{norm-diff}(a) the required estimate follows.
\end{proof}
\noindent By above Lemma and  the definition $p_t^{(h)}$ we have
\[ {\mathcal P}_t^{(h)}(x\textbf{e}(f))u= p_t^{(h)}(x) u\textbf{e}(f)= xu\textbf{e}(f)~~~~~~~~~~~~~~~~~~~~~~~~~~~~~~~~~~~~~\]
\be \label{qrw-mod}+\sum_{k=1}^n
  N_{p_{(k-1)h}^{(h)}({\beta(h,x)}- b(x))}[k] u\textbf{e}(f)+F(h,x,u,f)\ee
 \section{ EH flow as a strong limit of Quantum random walks} 
Here, we shall construct quantum random walk and  prove the  strong convergence  extending the ideas in \cite{KBS}.

Let $T_t$  be a  uniformly continuous conservative   QDS on von Neumann algebra   $\A$ with the generator ${\mathcal L}.$  Then (for detail see \cite{GS1}):\\
\noindent(i) There exists a Hilbert space  ${\mathbf{k}}_0$ and structure maps $({\mathcal L},\delta,  \sigma),$ where  $\delta \in \B(\A,\A \bigotimes {\mathbf{k}}_0)$ and $\sigma \in \B(\A,\A \bigotimes \B({\mathbf{k}}_0)).$\\
\noindent(ii) The map
 $\Theta=
\left( \begin{array}{cc} \theta_1 & (\theta_2(\cdot))^*\\
\theta_3 & \theta_4 \end{array}\right)
=\left( \begin{array}{cc} {\mathcal L} & \delta^\dag\\
\delta & \sigma \end{array}\right):\A \raro \A\bigotimes \B(\hat{{\mathbf{k}}_0})$ is a bounded CCP map
 (where $\delta^\dag(x)=(\delta(x^*))^*$ for all $x\in \A$)
 with the structure 
 \be \label{st-TH} \theta(x)=V^*(x\otimes 1_{\hat{\mathbf{k}}_0})V+ W(x\otimes 1_{\hat{\mathbf{k}}_0})+(x\otimes 1_{\hat{\mathbf{k}}_0})W^*,\forall x\in \A,\ee
where $V,W\in \B({\mathbf h}_0 \bigotimes \hat{\mathbf{k}}_0),$
 and the  estimate 
 \be \label{esti-TH}\|\Theta(x)\xi\|\le \|(x\otimes 1_{\H})D\xi\|, \forall x\in \A, \xi \in {\mathbf h}_0 \mbox{$ \bigotimes$} \hat{\mathbf{k}}_0,\forall x\in \A, \xi\in
  {\mathbf h}_0 \mbox{$ \bigotimes$} \hat{\mathbf{k}}_0,\ee 
  where $D\in \B({\mathbf h}_0 \bigotimes \hat{\mathbf{k}}_0,{\mathbf h}_0 \bigotimes \H),\H=\hat{\mathbf{k}}_0\bigoplus\hat{\mathbf{k}}_0\bigoplus \hat{\mathbf{k}}_0.$

\noindent (iii) Let $\tau \ge 0$ be fixed.  There exists a unique solution $J_t$ of the  equation,
 \be \label{Jt-exist}
 J_t=id_{\A \bigotimes \Gamma}+\int_0^t J_s\Lambda_{\Theta}(ds),~{0 \le t \le \tau} \ee
\[(\mbox{here we have written}~\Lambda_{\Theta}(ds) ~\mbox{for }~ \Lambda_{\theta_1}^1(ds)+\Lambda_{\theta_3}^2(ds)+\Lambda_{\theta_3}^3(ds)+\Lambda_{\theta_4}^4(ds))\]
as a regular adapted process mapping $\A \bigotimes \E(\C)$ into $\A \bigotimes \Gamma$ and satisfies
 $$ \sup_{0 \leq t \leq \tau} ||J_t(x \otimes \textbf{e}(f))u|| \leq
C^{\prime}(f) || (x \otimes 1_{\Gamma_{\rm fr}(L^2([0,\tau],\H))})E_{\tau}u||,$$
where $f \in \C, E_t \in  \B({\mathbf h}_0,{\mathbf h}_0\bigotimes  \Gamma_{\rm fr}(L^2([0,\tau],\H))),$
 $C^\prime(f)$ is some constant and \\$\Gamma_{\rm fr}(L^2([0,\tau],\H))$ is the free Fock space over $L^2([0,\tau],\H).$ \\

\noindent For $m\ge 0,$ let us consider the ampliation
\[\Theta_{(m)}:\A \mbox{$\bigotimes$} \B(\hat{\mathbf{k}}_0^{\mbox{\textcircled{\small m}}})\rightarrow \A 
\mbox{$\bigotimes$} \B(\hat{\mathbf{k}}_0^{\mbox{\textcircled{\small {m}}}})\mbox{$\bigotimes$} \B( \hat{\mathbf{k}}_0)~\mbox{of the map $\Theta$ given by}\]
\be \label{theta-m}\Theta_{(m)}(X)=Q_m^*\left(\Theta \otimes id_{\B(\hat{\mathbf{k}}_0^{\mbox{\textcircled{\small m}}})}(X)\right)Q_m\ee
 where $Q_m : {\mathbf h}_0 \bigotimes\hat{\mathbf{k}}_0^{\mbox{\textcircled{\small m}}}\bigotimes
 \hat{\mathbf{k}}_0\rightarrow {\mathbf h}_0 \bigotimes\hat{\mathbf{k}}_0
 \bigotimes  \hat{\mathbf{k}}_0^{\mbox{\textcircled{\small m}}}$ is  the unitary operator which interchanges the second and third tensor components.
  From the structure (\ref{st-TH}) of the map $\Theta,$ 
 \bean && \Theta_{(m)}(X)=Q_m^*(V^* \otimes 1_{\hat{\mathbf{k}}_0^{\mbox{\textcircled{\small m}}}}) Q_m (X\otimes 1_{\hat{\mathbf{k}}_0})Q_m^* (V\otimes 1_{\hat{\mathbf{k}}_0^{\mbox{\textcircled{\small m}}}})Q_m\\
 &&+Q_m^*( W\otimes 1_{\hat{\mathbf{k}}_0^{\mbox{\textcircled{\small m}}}}) Q_m(X\otimes 1_{\hat{\mathbf{k}}_0})+(X\otimes 1_{\hat{\mathbf{k}}_0})Q_m^*  (W^*\otimes 1_{\hat{\mathbf{k}}_0^{\mbox{\textcircled{\small m}}}}) Q_m.
 \eean
 For $\xi \in {\mathbf h}_0 \mbox{$\bigotimes$} \hat{\mathbf{k}}_0^{\mbox{\textcircled{\small m}}} \mbox{$\bigotimes$}\hat{\mathbf{k}}_0,$
 \bean
 &&\|\Theta_{(m)}(X)\xi\|^2
 \le 3\left [\|V\|^2 \|(X\otimes 1_{\hat{\mathbf{k}}_0})Q_m^* (V\otimes 1_{\hat{\mathbf{k}}_0^{\mbox{\textcircled{\small m}}}})Q_m \xi\|^2\right .\\
 &&\left .+\|W\|^2\|(X\otimes 1_{\hat{\mathbf{k}}_0})\xi \|^2 
 +\|(X\otimes 1_{\hat{\mathbf{k}}_0})Q_m^*  (W^*\otimes 1_{\hat{\mathbf{k}}_0^{\mbox{\textcircled{\small m}}}}) Q_m \xi\|^2\right ].
 \eean
 Setting \[D_m\xi=\sqrt{3}\left [\|V\| Q_m^*(V\otimes 1_{\hat{\mathbf{k}}_0^{\mbox{\textcircled{\small m}}}})Q_m \xi
 \oplus \|W\|\xi 
 \oplus Q_m^*  (W^*\otimes 1_{\hat{\mathbf{k}}_0^{\mbox{\textcircled{\small m}}}}) Q_m \xi\right ],\]
 $D_m\in \B({\mathbf h}_0 \mbox{$\bigotimes$} \hat{\mathbf{k}}_0^{\mbox{\textcircled{\small m}}} \mbox{$\bigotimes$}\hat{\mathbf{k}}_0,{\mathbf h}_0 \mbox{$\bigotimes$} \hat{\mathbf{k}}_0^{\mbox{\textcircled{\small m}}} \mbox{$\bigotimes$}\H)$
 (where $\H=\hat{\mathbf{k}}_0\oplus \hat{\mathbf{k}}_0\oplus \hat{\mathbf{k}}_0$ as earlier) and
 \be \label{Theta-ampli}\|\Theta_{(m)}(X)\xi\|\le \|(X\otimes 1_{\H})D_m\xi\|, \forall X\in \A \mbox{$\bigotimes$} \B(\hat{\mathbf{k}}_0^{\mbox{\textcircled{\small m}}}).\ee
 Thus $\|\Theta_{(m)}\|\le \|D_m\|,$
by definition $\|D_m\|^2\le 3(\|V\|^4+\|W\|^2), \forall m\ge 0$ and hence $\Theta $ can be  extend as a map 
 $\bigoplus_{m\ge 0} \Theta_{(m)}$ from $ \A \bigotimes \B(\Gamma_{\rm fr}(\hat{\mathbf{k}}_0))$ into itself  with
$\|\bigoplus_{m\ge 0}  \Theta_{(m)}\|\le 3(\|V\|^2+\|W\|),$ we denote this map by  same symbol $\Theta.$

\noindent For any fixed  $m\ge 0$ let us look at the following  qsde on $\A \bigotimes \B(\hat{\mathbf{k}}_0^{\mbox{\textcircled{\small m}}})
\bigotimes \Gamma $
\be \label{mod-qsde}
\eta_{m,t} =id_{\A \bigotimes \B(\hat{\mathbf{k}}_0^{\mbox{\textcircled{\small m}}})
\bigotimes \Gamma }+\int_0^t \eta_{m,s}  \Lambda_{\Theta}(ds), 0\le t\le \tau.
\ee
Since we have the estimate, for any $ X\in \A \bigotimes \B(\hat{\mathbf{k}}_0^{\mbox{\textcircled{\small m}}}),~ \xi \in {\mathbf h}_0 \bigotimes \hat{\mathbf{k}}_0^{\mbox{\textcircled{\small m}}} \bigotimes \hat{\mathbf{k}}_0$
\[\|\Theta(X)\xi\|=\|\Theta _{(m)}(X)\xi\|\le \|(X\otimes 1_{\H})D_m\xi\|, \] by a simple adoptation of the proof of the existence of solution $J_t$ of the qsde \ref{Jt-exist} (Theorem 3.3.6 (i) in \cite{GS1}), it can be shown that\\ \noindent(i) the qsde (\ref {mod-qsde}) admit a unique solution
 $\eta_{m,t}$  as an adapted regular process mapping $\A \bigotimes \B(\hat{\mathbf{k}}_0^{\mbox{\textcircled{\small m}}})
\bigotimes \E(\C)$ into $\A \bigotimes \B(\hat{\mathbf{k}}_0^{\mbox{\textcircled{\small m}}})
\bigotimes \Gamma.$ \\
\noindent(ii) $\eta_{m,t}$ satisfies the  
estimate
 \be \label{mod-esti} \sup_{0 \leq t \leq \tau} ||\eta_{m,t}(X \otimes \textbf{e}(f))\xi|| \leq
C^{\prime}(f) || (X \otimes 1_{\Gamma_{\rm fr}(L^2([0,\tau],\H))})E_{\tau}\xi||,\ee
where $f \in \C$ and 
 $C^\prime(f)$ is some constant. The operator  $E_{\tau}$ appears above is  an element of \\
 $ \B({\mathbf h}_0\bigotimes \hat{\mathbf{k}}_0^{\mbox{\textcircled{\small m}}},{\mathbf h}_0\bigotimes \hat{\mathbf{k}}_0^{\mbox{\textcircled{\small m}}}\bigotimes  \Gamma_{\rm fr}(L^2([0,\tau],\H))),$ define as follows:\\
 \[ E_{\tau}\xi= \bigoplus_{n\ge 0} (n!)^{\frac{1}{4}} E_{\tau}^{(n)}\xi,\]
 where $ E_{\tau}^{(n)}\in \B({\mathbf h}_0\bigotimes \hat{\mathbf{k}}_0^{\mbox{\textcircled{\small m}}},{\mathbf h}_0\bigotimes \hat{\mathbf{k}}_0^{\mbox{\textcircled{\small m}}}\bigotimes  (L^2([0,\tau],\H))^{\mbox{\textcircled{\small n}}})$ given by, for $\xi\in {\mathbf h}_0 \bigotimes \hat{\mathbf{k}}_0^{\mbox{\textcircled{\small m}}}$
 \bean
 && E_{\tau}^{(0)}\xi =\xi,\\
 && (E_{\tau}^{(1)}\xi )(s) =
   D(\xi \otimes\hat{f}(s)||\hat{f_{\left. t \right]}}(s)||)
  ~\mbox{ and iteratively}\\ 
 &&(E_{\tau}^{(n)}\xi)(s_1,s_2,\ldots s_n)=(D_m \otimes 1_{{L^2([0,\tau],\H)}^{\otimes^{n-1}}}){\mathcal Q}_n \\
 &&~~~~~~~~~~~~~~~~~~~~~~~~~~~~~~~\{ (E_{\tau}^{(n-1)}u)(s_2, \ldots s_n)\otimes\hat{f}(s_1)||\hat{f_{\left. t \right]}}(s_1)||\} 
 \eean
 ( ${\mathcal Q}_n : {\mathbf h}_0 \bigotimes \hat{\mathbf{k}}_0^{\mbox{\textcircled{\small m}}} \bigotimes{L^2([0,\tau],\H)}^{\otimes^{(n-1)}} \bigotimes
 \hat{\mathbf k}_0 \rightarrow {\mathbf h}_0 \bigotimes \hat{\mathbf{k}}_0^{\mbox{\textcircled{\small m}}} \bigotimes 
 \hat{\mathbf k}_0\bigotimes{L^2([0,\tau],\H)}^{\otimes^{(n-1)}}$ is  the unitary operator which interchanges the  third and fourth tensor components).
 
  It is clear that $J_t\otimes id_{\B(\hat{\mathbf{k}}_0^{\mbox{\textcircled{\small m}}})}
   ( \equiv \Upsilon_m^*(J_t\otimes id_{\B(\hat{\mathbf{k}}_0^{\mbox{\textcircled{\small m}}})})\Upsilon_m:\A \bigotimes \B(\hat{\mathbf{k}}_0^{\mbox{\textcircled{\small m}}})\bigotimes \E(\C)\raro
\A \bigotimes \B(\hat{\mathbf{k}}_0^{\mbox{\textcircled{\small m}}})\bigotimes \Gamma,$ where $\Upsilon_m:{\mathbf h}_0 \bigotimes \hat{\mathbf{k}}_0^{\mbox{\textcircled{\small m}}} \bigotimes \Gamma\raro
{\mathbf h}_0 \bigotimes \Gamma \bigotimes \hat{\mathbf{k}}_0^{\mbox{\textcircled{\small m}}}$)\\
 satisfies the
qsde (\ref {mod-qsde}) and hence $\eta_{m,t}=J_t\otimes id_{\B(\hat{\mathbf{k}}_0^{\mbox{\textcircled{\small m}}})}.$ By definition of $E_{\tau},$ it can be easily seen that $\|E_{\tau}\|$ uniformly bounded for $m\ge 0 $ and hence  the estimate (\ref{mod-esti}) allow us to 
extend $\{J_t\}$ as a regular adapted process  $\{\bigoplus_{m\ge 0} J_t\otimes id_{\B(\hat{\mathbf{k}}_0^{\mbox{\textcircled{\small m}}})}\}$ mapping  $\A \bigotimes \B(\Gamma_{\rm fr}(\hat{\mathbf{k}}_0))\bigotimes \E(\C)$ into
$\A \bigotimes \B(\Gamma_{\rm fr}(\hat{\mathbf{k}}_0))\bigotimes \Gamma,$ we denote this family by same symbol $J_t.$ For a given $f\in \C$  this $J_t$ satisfies
 \be \label{bigJtesti}
 \|J_t(X\textbf{e}(f)\xi\|\le D^\prime \|X\| \|\xi\|, \forall X\in \A \bigotimes \B(\Gamma_{\rm fr}(\hat{\mathbf{k}}_0)) ~\mbox{and}~ \xi\in {\mathbf h}_0 \bigotimes \Gamma_{\rm fr}(\hat{\mathbf{k}}_0)), 
 \ee
 for some constant $D^\prime$ independent of $X$ and $\xi.$

\section*{Construction of quantum random walks}
In particular suppose we are given with a conservative  uniformly continuous QDS $T_t$ on a von Neumann algebra $\A$ with generator ${\mathcal L}$  given by
\be \label{LL1}{\mathcal L}(x)=R^* (x\otimes 1_{{\mathbf{k}}_0}) R-\half R^*Rx -\half xR^*R, \forall x\in \A  \ee for some Hilbert space ${\mathbf{k}}_0$ and $R\in \A \bigotimes {\mathbf{k}}_0.$ 

\bthm \label{Beta-Tt}
Let ${\mathcal L}$ be given by (\ref{LL1}). Then there exists a $*$-homomorphic family $\{\beta(h)\}_{h> 0}$ such that  the family of linear   maps  $E(h):\A\raro \A\bigotimes \B(\hat{\mathbf{k}}_0)$ given by, for $x\in \A$  \\
$ E(h,x)
=\left( \begin{array}{cc}  h^{-2}\left [\beta_1(h,x)- x-h
 \theta_1(x) \right ] &  h^{-\frac{3}{2}} (\beta_2(h,x)-\sqrt{h}
 \theta_2(x))^*\\
h^{-\frac{3}{2}} [\beta_3(h,x)-\sqrt{h}
 \theta_3(x) ] & h^{-1} \left [\beta_4(h,x)- x\otimes 1_{{\mathbf{k}}_0}-
 \theta_4(x) \right ]\end{array}\right),$\\ is  uniformly norm bounded, also  as   maps from 
 $\A \bigotimes \B(\Gamma_{\rm{fr}}(\hat{\mathbf{k}}_0))$ into itself, have  uniform norm bound
 i.e. $\|E(h)\| \le M,$ for some constant $M$ independent of $h.$ 
\ethm   
\noindent In particular, it follows  that for any $l$ 
  \be \label{theta-sig}
 \|\beta_l(h,X)- b_l(X)-h^{\varepsilon_l}
 \theta_l(X)\|
 \le M \|X\| h^{1+\varepsilon_l}, \forall X\in \A ~\mbox{$\bigotimes$}\B(\hat{\mathbf{k}}_0^{\mbox{\textcircled{\small m}}})\ee
 where $\varepsilon_1=1,\varepsilon_2=\varepsilon_3=\frac{1}{2}$ and
$\varepsilon_4=0.$

\begin{proof}
Here  the map  $\Theta $ is  given by  $\Theta(x) =\left( \begin{array}{cc} \theta_1(x) & (\theta_2(x))^*\\
\theta_3(x) &  \theta_4(x)
\end{array} \right)=\left( \begin{array}{cc} {\mathcal L} (x) & \delta^\dag(x)\\
\delta(x) &  \sigma(x)
\end{array} \right),\\  \forall x\in \A,$ where $\delta(x)= (x\otimes 1_{{\mathbf{k}}_0}) R-R
x ,\delta^\dag(x)=(\delta(x^*))^*=R^*( x\otimes 1_{{\mathbf{k}}_0})-xR^*$  and $\sigma=0.$ 
Setting $\widetilde{R}=
\left( \begin{array}{cc} 0 & -R^*\\
 R & 0
\end{array} \right)$ from ${\mathbf{h}}_0 \bigotimes \hat{\mathbf{k}}_0$ to itself.  It is clear that $ \widetilde{R}$ is a  bounded skew symmetric  operator thus it generate a one parameter  unitary group $\{e^{t \widetilde{R} }\}.$ For $h>0,$ we consider the unitary operator $U(h)=e^{\sqrt{h} \widetilde{R} }$
which can be written as
$\left( \begin{array}{cc} \cos(\sqrt{h}|R|)&
-\sqrt{h} D(h)R^*\\
\sqrt{h}R D(h) & \cos(\sqrt{h}|R^*|)
\end{array} \right)$ where 
$D(h)=\sin(\sqrt{h} |R|){( \sqrt{h}|R|)}^{-1}$ and  $|R|,|R^*|$
 denote the  positive square root of $R^*R$ and $RR^*$ respectively.
It can easily be   observed that 
 \begin{equation} \label{un-esti}
 \begin{split}
&\|\cos(\sqrt{h}|R|)- 1_{{\mathbf h}_0} + \frac{h}{2} |R|^2\| \le h^2 \|R\|^4,\\
&\|\cos(\sqrt{h}|R|)-1_{{\mathbf h}_0 }\| \le h \|R\|^2,\\
&\|\cos(\sqrt{h}|R^*|)-1_{{\mathbf h}_0\otimes {\mathbf{k}}_0}\|\le h \|R\|^2,\\
&\|D(h)-1_{{\mathbf h}_0}\| \le h \|R\|^2,\\
&\|\cos(\sqrt{h}|R|)\|\le 1,\\
&\|D(h)\| \le 1.
\end{split}
\end{equation}
Now we define a $*$-homomorphism  $\beta(h)$ from $\A$ to $\A\bigotimes
\B(\hat{\mathbf{k}}_0)$ implemented by the unitary $U(h),$ i.e.
for $x\in \A,\beta(h,x):=\beta(h)(x)=(U(h))^*
(x \otimes 1_{\hat{\mathbf{k}}_0})U(h).$  So   for any  $x\in \A,
\beta(h,x) =\left( \begin{array}{cc} \beta_1(h,x) & (\beta_2(h,x))^*\\
\beta_3(h,x) &  \beta_4(h,x)
\end{array} \right)
\\
=\left( \begin{array}{cc}
\{\cos(\sqrt{h}|R|)x \cos(\sqrt{h}|R|)
&\{-\sqrt{h}\cos(\sqrt{h}|R|)x D(h)R^*\\
\vspace{5mm}+h D(h)R^* (x\otimes 1_{{\mathbf{k}}_0})R D(h)\}&
+\sqrt{h} D(h)R^*(x\otimes 1_{{\mathbf{k}}_0})\cos(\sqrt{h}|R^*|)\}\\
 \{-\sqrt{h}R D(h) x\cos(\sqrt{h}|R|)& \{h R
D(h) x D(h)R^*\\
+\sqrt{h}\cos(\sqrt{h}|R^*|)(x\otimes 1_{{\mathbf{k}}_0}) R D(h)\}& +
\cos(\sqrt{h}|R^*|) (x\otimes 1_{{\mathbf{k}}_0}) \cos(\sqrt{h}|R^*|)\}
\end{array} \right).$
We have
\bean
&&\beta_1(h,x)-x-h \theta_1(x)\\
&&=\cos(\sqrt{h}|R|)x \cos(\sqrt{h}|R|)+h D(h)R^* (x\otimes 1_{{\mathbf{k}}_0})R D(h)\\
&&~~-x-h\left(R^* (x\otimes 1_{{\mathbf{k}}_0}) R-\half |R|^2x -\half x|R|^2\right)\\
&&=\left [\cos(\sqrt{h}|R|)-1_{{\mathbf h}_0}+\half |R|^2 \right ]x \cos(\sqrt{h}|R|)\\
&&~~+x\left [\cos(\sqrt{h}|R|)-1_{{\mathbf h}_0}+\half |R|^2 \right ]
+\half |R|^2x \left [1_{{\mathbf h}_0}-\cos(\sqrt{h}|R|)\right ]\\
&&~~+h[ D(h)-1_{{\mathbf h}_0}]R^* (x\otimes 1_{{\mathbf{k}}_0})R D(h)
+hR^* (x\otimes 1_{{\mathbf{k}}_0})R D(h).
\eean
By (\ref{un-esti}) we get
\be \label{esti1}
\|\beta_1(h,x)-x-h \theta_1(x)\|\le 5  h^2 \|R\|^4 \|x\|.\\
\ee
By definition we have 
\bean
&& \beta_2(x^*)-\sqrt{h}\theta_2(x^*)=\beta_3(x)-\sqrt{h}\theta_3(x)\\
&&= \sqrt{h}\left [-R D(h) x\cos(\sqrt{h}|R|)+\cos(\sqrt{h}|R^*|)(x\otimes 1_{{\mathbf{k}}_0}) R D(h) - (x\otimes 1_{{\mathbf{k}}_0}) R+ R x\right ]\\
&&=\sqrt{h} \left [-R D(h) x[\cos(\sqrt{h}|R|)-1_{{\mathbf h}_0}]-R [D(h)-1_{{\mathbf h}_0}] x\right.\\
&&~~~~+\left .\cos(\sqrt{h}|R^*|)(x\otimes 1_{{\mathbf{k}}_0}) R [D(h)-1_{{\mathbf h}_0}]
+ [\cos(\sqrt{h}|R^*|)-1_{{\mathbf h}_0}](x\otimes 1_{{\mathbf{k}}_0}) R\right ].
\eean
Using (\ref{un-esti})  we get
\bean
&& \|\beta_2(x^*)-\sqrt{h}\theta_2(x^*)\|=\|\beta_3(x)-\sqrt{h}\theta_3(x)\|\\
&&\le \sqrt{h}\left[ \|R D(h) x[\cos(\sqrt{h}|R|)-1_{{\mathbf h}_0}]\|
+\|R [D(h)-1_{{\mathbf h}_0}] x\|\right.\\
&&~~~~\left.+\|\cos(\sqrt{h}|R^*|)(x\otimes 1_{{\mathbf{k}}_0}) R [D(h)-1_{{\mathbf h}_0}]\|
+\| [\cos(\sqrt{h}|R^*|)-1_{{\mathbf h}_0\otimes {\mathbf{k}}_0}](x\otimes 1_{{\mathbf{k}}_0}) R\|\right]
\\
&&~~~\le 4 h^{\frac{3}{2}} \|R\|^3 \|x\|.
\eean
Now let us consider $\beta_4(x)-\theta_4(x),$ we have
\bean 
&& \|\beta_4(x)-\theta_4(x)\|\\
&&= \|h RD(h) x D(h)R^* +
\cos(\sqrt{h}|R^*|) (x\otimes 1_{{\mathbf{k}}_0}) \cos(\sqrt{h}|R^*|)-(x\otimes 1_{{\mathbf{k}}_0})\|\\
&& \le   h\| RD(h) x D(h)R^* \|+
\|[\cos(\sqrt{h}|R^*|)-1_{{\mathbf h}_0\otimes {\mathbf{k}}_0}] (x\otimes 1_{{\mathbf{k}}_0}) \cos(\sqrt{h}|R^*|)\| \\
&&~~+\|(x\otimes 1_{{\mathbf{k}}_0}) [\cos(\sqrt{h}|R^*|)-1_{{\mathbf h}_0\otimes {\mathbf{k}}_0}]\|\\
&&\le 3 h \|R\|^2 \|x\|.
\eean
Thus for $l=1,2,3$ and $4,$ 
\be \label{thetaesti}
\|\beta_l(h,x)- b_l(x)-h^{\varepsilon_l}
 \theta_l(x)\|
 \le M \|x\| h^{1+\varepsilon_l},  \forall  x\in \A,
 \ee
 where constant $M=5(\|R\|^2+\|R\|^3+\|R\|^4).$

 \noindent For $m\ge 0,$ let us consider the ampliation of the maps  $\Theta,b$ and $\beta$ as maps from $\A \bigotimes  \B(\Gamma_{\rm fr}(\hat{\mathbf{k}}_0))$
into itself. For $X\in \A \mbox{$\bigotimes$} \B(\hat{\mathbf{k}}_0^{\mbox{\textcircled{\small m}}})$
\[ \Theta (X)=\Theta_{(m)}(X)=Q_m^*\left(\Theta \otimes id_{\B(\hat{\mathbf{k}}_0^{\mbox{\textcircled{\small m}}})}(X)\right)Q_m,\]
 where $Q_m : {\mathbf h}_0 \bigotimes\hat{\mathbf{k}}_0^{\mbox{\textcircled{\small m}}}\bigotimes
 \hat{\mathbf{k}}_0\rightarrow {\mathbf h}_0 \bigotimes\hat{\mathbf{k}}_0
 \bigotimes  \hat{\mathbf{k}}_0^{\mbox{\textcircled{\small m}}}$ is  the unitary operator which interchanges the second and third tensor components. This operator 
 $Q_m=1_{{\mathbf h}_0\bigotimes \hat{\mathbf{k}}_0^{\mbox{\textcircled{\small m}}}} \oplus P_m$ where
 $q_m : {\mathbf h}_0 \bigotimes\hat{\mathbf{k}}_0^{\mbox{\textcircled{\small m}}}\bigotimes
 {\mathbf{k}}_0\rightarrow {\mathbf h}_0 \bigotimes{\mathbf{k}}_0
 \bigotimes  \hat{\mathbf{k}}_0^{\mbox{\textcircled{\small m}}}$ is define as $Q_m.$ By definition we have
\bean
&&\theta_1(X)=(R^* \otimes 1_{\hat{\mathbf{k}}_0^{\mbox{\textcircled{\small m}}}}) q_m(X\otimes 1_{{\mathbf{k}}_0})q_m^* (R\otimes 1_{\hat{\mathbf{k}}_0^{\mbox{\textcircled{\small m}}}})-\half (|R|\otimes 1_{\hat{\mathbf{k}}_0^{\mbox{\textcircled{\small m}}}})X -\half X(|R|\otimes 1_{\hat{\mathbf{k}}_0^{\mbox{\textcircled{\small m}}}})\\
&&
\theta_2(X)^*=\left[ X(R^* \otimes 1_{\hat{\mathbf{k}}_0^{\mbox{\textcircled{\small m}}}}) q_m-(R^* \otimes 1_{\hat{\mathbf{k}}_0^{\mbox{\textcircled{\small m}}}}) q_m(X\otimes 1_{{\mathbf{k}}_0})\right ]q_m^*\\
&&
\theta_3(X)=q_m^*\left[(R \otimes 1_{\hat{\mathbf{k}}_0^{\mbox{\textcircled{\small m}}}})X- q_m(X\otimes 1_{{\mathbf{k}}_0})q_m^* (R\otimes 1_{\hat{\mathbf{k}}_0^{\mbox{\textcircled{\small m}}}})\right]\\
&&
\theta_4(X)=0
\eean

and components of  $\beta(h,X)$ are

\bean
&&
\beta_1(h,X)=(\cos(\sqrt{h}|R|) \otimes 1_{\hat{\mathbf{k}}_0^{\mbox{\textcircled{\small m}}}})X (\cos(\sqrt{h}|R|)\otimes 1_{\hat{\mathbf{k}}_0^{\mbox{\textcircled{\small m}}}})\\
&&~~~~~~~~
+h (D(h)R^*\otimes 1_{\hat{\mathbf{k}}_0^{\mbox{\textcircled{\small m}}}}) q_m(X\otimes 1_{{\mathbf{k}}_0})q_m^*(R D(h)\otimes 1_{\hat{\mathbf{k}}_0^{\mbox{\textcircled{\small m}}}})\\
&&
\beta_2(h,X)^*=\left [-\sqrt{h}(\cos(\sqrt{h}|R|)\otimes 1_{\hat{\mathbf{k}}_0^{\mbox{\textcircled{\small m}}}})X (D(h)R^*\otimes 1_{\hat{\mathbf{k}}_0^{\mbox{\textcircled{\small m}}}})\right.\\
&&~~~~~~~~
\left. +\sqrt{h} (D(h)R^*\otimes 1_{\hat{\mathbf{k}}_0^{\mbox{\textcircled{\small m}}}})q_m(X\otimes 1_{{\mathbf{k}}_0})q_m^*(\cos(\sqrt{h}|R^*|)\otimes 1_{\hat{\mathbf{k}}_0^{\mbox{\textcircled{\small m}}}})\right ] q_m\\
&&
\beta_3(h,X)=q_m^*\left [-\sqrt{h}(R D(h)\otimes 1_{\hat{\mathbf{k}}_0^{\mbox{\textcircled{\small m}}}}) X(\cos(\sqrt{h}|R|)\otimes 1_{\hat{\mathbf{k}}_0^{\mbox{\textcircled{\small m}}}})\right. \\
&&~~~~~~~~
\left. +\sqrt{h}(\cos(\sqrt{h}|R^*|)\otimes 1_{\hat{\mathbf{k}}_0^{\mbox{\textcircled{\small m}}}})q_m(X\otimes 1_{{\mathbf{k}}_0})q_m^* (R D(h)\otimes 1_{\hat{\mathbf{k}}_0^{\mbox{\textcircled{\small m}}}})\right]\\
&&
\beta_4(h,X)=q_m^*\left [h (R
D(h)\otimes 1_{\hat{\mathbf{k}}_0^{\mbox{\textcircled{\small m}}}}) X ( D(h)R^*
\otimes 1_{\hat{\mathbf{k}}_0^{\mbox{\textcircled{\small m}}}})\right.\\
&&~~~~~~~~
\left . +(\cos(\sqrt{h}|R^*|)\otimes 1_{\hat{\mathbf{k}}_0^{\mbox{\textcircled{\small m}}}}) q_m(X\otimes 1_{{\mathbf{k}}_0})q_m^* (\cos(\sqrt{h}|R^*|)\otimes 1_{\hat{\mathbf{k}}_0^{\mbox{\textcircled{\small m}}}})\right ] q_m.
\eean
\noindent By same argument  as
for (\ref{thetaesti})  one has
\[\|\beta_1 (h,X)- b_l(X) -h \theta_1 (X)\|\le C^\prime h^{1+\varepsilon_l}  \|X\|, \forall  X\in \A \bigotimes \B(\hat{\mathbf{k}}_0^{\mbox{\textcircled{\small m}}}),\] for some constant $C^\prime$ independent of $h >0 , m\ge 0.$ 
Thus $\|E(h)(X)\|\le M \|X\|, \forall X\in \A\bigotimes \B(\Gamma_{\rm fr}(\hat{\mathbf{k}}_0)),$ for some constant $M$ independent  of $h.$  
\end{proof}

Now consider the  quantum random walks $\{p_t^{(h)}:h> 0\}$ associated with 
the $*$-homomorphic family $\{\beta(h)\}$ in the above theorem. In the next section we shall prove that this  quantum random walks $\{p_t^{(h)}:h> 0\}$
converges strongly.

\section*{Convergence of quantum random walk}
\noindent Let $\B=\A\bigotimes \B(\Gamma_{\rm{fr}}(\hat{\mathbf{k}}_0))\bigotimes \B(\Gamma),$ which can be decomposed as \\
$\B=\left (\A\bigotimes  \bigoplus_{m\ge 0} \B(\hat{\mathbf{k}}_0^{\mbox{\textcircled{\small m}}})\right ) \bigoplus \B_c$ for some subspace $\B_c.$ Now let us consider the extensions of  all these  maps $\Theta ,\beta(h), b,p_t^{(h)}$
 and ${\mathcal P}_t^{(h)}$ as bounded linear  maps from $\B$ into itself, given by, for example extention of $p_t^{(h)}$ is
  $\bigoplus_{m\ge 0}~ p_t^{(h)} \otimes id_{\B(\hat{\mathbf{k}}_0^{\mbox{\textcircled{\small{m}}}})} \oplus 0_{\B_c}.$
 We denote these extentions by  same symbols as the original maps.
 \noindent From the Theorem \ref{Beta-Tt} it follows that 
\be \label{theta-cb} 
 \|\beta_l(h,X)- b_l(X)-h^{\varepsilon_l}
 \theta_l(X)\|
 \le C \|X\| h^{1+\varepsilon_l},  \forall  X\in \A \mbox{$\bigotimes$} \B(\hat{\mathbf{k}}_0^{\mbox{\textcircled{\small m}}}).\ee
For $h\ge 0$ we define a map
$\Theta(h):=
\left( \begin{array}{cc} h\theta_1 & h^\half(\theta_2(\cdot))^*\\
h^\half\theta_3 & \theta_4 \end{array}\right)$ from $\A$ to $\A\bigotimes \B(\hat{\mathbf{k}}_0),$ as the map $\Theta, ~\Theta(h)$ also extend as a bounded map from $\B$ into itself.
Here we  have the  following observations which will be needed later for proving the  convergence of quantum random walk $p_{t}^{(h)}.$
\blema \label{NLamda} For any $l, X\in\A\bigotimes \B(\hat{\mathbf{k}}_0^{\mbox{\textcircled{\small m}}}), \xi \in {\mathbf h}_0\bigotimes
\hat{\mathbf{k}}_0^{\mbox{\textcircled{\small m}}}$ and $ f\in \M$ we have
\begin{enumerate}
\item $  \| \sum_{k=1}^n 
  N_{p_{(k-1)h}^{(h)} [\beta(h,X)- b(X)
-\Theta(h,X)]}   [k]\xi \textbf{e}(f)\| \le \sqrt{h} C_1(f,t)\|X\| \|\xi\|$\\
\item $ \|\sum_{k=1}^n
 \left [ N_{p_{(k-1)h}^{(h)}(\Theta(h,X))}[k]-  \Lambda_{p_{(k-1)h}^{(h)} (\Theta(X))} [k] \right ]\xi \textbf{e}(f)\|^2\\ \\ \le h C_2(f,t)\|X\|^2 \|\xi \|^2, $
 \item $\|\sum_{k=1}^n p_{(k-1)h}^{(h)}(X)(1-P_h[k])\xi\textbf{e}(f))\|^2
 \le h~c(f,t)\|X\|^2 \|\xi\|^2.$
\end{enumerate}
where  constants $c(f,t)$ is as in Lemma \ref{JPh}, $C_1(f,t)=t (1+\|f\|_\infty)\|\textbf{e}(f)\|$
and\\
$C_2(f,t)=(1+t)^2 (\|f\|_\infty +\|f\|_\infty^2)^2 (1+\|\Theta\|)^2 \|\textbf{e}(f)\|^2.$
\elema       
 \begin{proof}
 (1). For any $l$ we have 
\bean &&  \|\sum_{k=1}^n  N_{p_{(k-1)h}^{(h)} [\beta_l(h,X)- b_l(X)
-h^{\varepsilon_l}\theta_l(X)]} ^l  [k]\xi \textbf{e}(f)\|\\
&& \le  \sum_{k=1}^n 
  \| N_{p_{(k-1)h}^{(h)} [\beta_l(h,X)- b_l(X)
-h^{\varepsilon_l}\theta_l(X)]} ^l  [k]\xi_{k-1} \textbf{e}(f_{[k]})\| \|\textbf{e}(f_{\left[kh\right .})\|
\eean
where  $\xi_{k-1}=\xi\textbf{e}(f_{(k-1)h]})$ is a  vector in the initial Hilbert space \\ ${\mathbf h}_0 \bigotimes \Gamma_{\rm{fr}}(\hat{\mathbf{k}}_0)\bigotimes \Gamma_{\left .(k-1)h\right ]}.$  For any $l,$   from (\ref{theta-cb}) and contractivity of $p_t^{(h)},$ we get
\[\|p_{(k-1)h}^{(h)} [\beta_l(h,X)- b_l(X)
-h^{\varepsilon_l}\theta_l(X)]\|\le C h^{1+\varepsilon_l} \|X\|,\]  hence by  ~(\ref{N-norm}) the above quantity is dominated by\\ $  \sum_{k=1}^n h^{\frac{3}{2}} C (1+\|f\|_\infty) \|X\|~ \|\xi \textbf{e}(f)\|$ and required estimate follows.\\
\noindent (2). By Lemma ~\ref{N-lambda}~the terms  correspond to  $l=1,2$
can be estimated as,
\bean
\lefteqn{\|\sum_{k=1}^n
\left [h^{\varepsilon_l} N_{p_{(k-1)h}^{(h)} (\theta_l(X))} ^l[k]-  \Lambda_{p_{(k-1)h}^{(h)}(\theta_l(X)) }^l[k]\right ]\xi \textbf{e}(f)\|}\\
&\le & \sum_{k=1}^n
 \|\left [h^{\varepsilon_l} N_{p_{(k-1)h}^{(h)}(\theta_l(X))} ^l[k]-  \Lambda_{p_{(k-1)h}^{(h)}(\theta_l(X))} ^l[k]\right ]\xi_{k-1} \textbf{e}(f_{[k]})\|
 \|\textbf{e}(f_{\left [kh \right .})\|\\
 &\le & \sum_{k=1}^n h^{\frac{3}{2}}
 \|p_{(k-1)h}^{(h)}(\theta_l(X))\| \|\xi_{k-1}\| (\|f\|_\infty + \|f\|_\infty^2 ) \|\textbf{e}(f_{[(k-1)h})\|~\\
 &\le & (\|f\|_\infty + \|f\|_\infty^2 ) \|\Theta\|\sum_{k=1}^n h^{\frac{3}{2}}
 \|X\| \|\xi \textbf{e}(f)\|.
\eean
Thus the required estimate follows.
Now consider other two  terms  correspond  to $l=3$ and $4.$ 
Setting for $1\le m\le n$
\[Z_{m}=\sum_{k=1}^m
 \left [\sqrt{h} N_{p_{(k-1)h}^{(h)}(\theta_l(X))}^l[k]-  \Lambda_{p_{(k-1)h}^{(h)}(\theta_l(X))} ^l[k]\right ],\]
by Lemma~ \ref{N-lambda} (a), we have
 \bean \lefteqn{\|Z_{m} u\textbf{e}(f_{mh]}) \|}\\
 &\le & \sum_{k=1}^m
 \|\left [\sqrt{h} N_{p_{(k-1)h}^{(h)}(\theta_l(X))} ^l[k]-  \Lambda_{p_{(k-1)h}^{(h)}(\theta_l(X))} ^l[k]\right ] \xi_{(k-1)}\textbf{e}(f_{[k]})\| \|\textbf{e}(f_{(kh, mh]})\|\\
 &\le & \sum_{k=1}^m h \|p_{(k-1)h}^{(h)}(\theta_l(X))\xi_{k-1}  \|
~\|f\|_\infty ~\|\textbf{e}(f_{[k]})\| \|\textbf{e}(f_{(kh, mh]})\| .
\eean
Thus
\be \label{ZZ}\|Z_{m} u\textbf{e}(f_{mh]}) \|\le  t \|\Theta\| ~\|f\|_\infty \|X\| \|\xi  \textbf{e}(f_{mh]}) \|.\ee
 We have the following equality,
 \bean \lefteqn{\|Z_n \xi \textbf{e}(f)\|^2}\\
 &=&\sum_{k=1}^n\|\left [\sqrt{h} N_{{(k-1)h}^{(h)}(\theta_l(X))} ^l[k]-  \Lambda_{p_{(k-1)h}^{(h)}(\theta_l(X))} ^l[k]\right ]\xi_{(k-1)} \textbf{e}(f_{[k]})\|^2  \textbf{e}(f_{\left [kh\right .})\|^2\\
 && +2 ~\mathcal {R}e~ \sum_{k=1}^n\lgl Z_{k-1} \xi_{k-1} \textbf{e}(f_{[k]})  ,~   \left [\sqrt{h} N_{p_{(k-1)h}^{(h)}(\theta_l(X))} ^l[k]-  \Lambda_{p_{(k-1)h}^{(h)}(\theta_l(X))}^l[k]\right ]\\
 &&~~~~~~~~~~~~~~~~~~~~~~\xi_{k-1} \textbf{e}(f_{[k]})  \rgl
 \|\textbf{e}(f_{\left[kh\right .})\|^2.
 \eean
By the estimate in Lemma~ \ref{N-lambda}, 
\bean
\lefteqn{\|\sum_{k=1}^n
 \left [\sqrt{h} N_{p_{(k-1)h}^{(h)}(\theta_l(X))} ^l[k]-  \Lambda_{p_{(k-1)h}^{(h)}(\theta_l(X))}^l[k]\right ]\xi \textbf{e}(f)\|^2}\\
&\le &\sum_{k=1}^n h^2 \|p_{(k-1)h}^{(h)}(\theta_l(X))\xi_{k-1}\|^2
 \|f\|_\infty^2  \|\textbf{e}(f_{[(k-1)h})\|^2\\
 &&+2 \sum_{k=1}^n h^2  \|Z_{k-1} \xi_{k-1}\|   \|p_{(k-1)h}^{(h)}(\theta_l(X))\|~\| \xi_{k-1} \|
 ~\|f\|_{\infty}^2 \|\textbf{e}(f_{[(k-1)h})\|^2.
 \eean
 Now using (\ref{ZZ}), above quantity is less than or equal to
 \bean &&\sum_{k=1}^n h^2  \|f\|_\infty^2~\|X\|^2 \|\xi \textbf{e}(f)\|^2  \\
 &&+2 \sum_{k=1}^n h^2  t \|\Theta\|^2  ~\|f\|_{\infty}^3 \|X\|^2 
 \|\xi \textbf{e}(f)\|^2 
 \eean
 and required estimate follows.\\
 (3.) The proof is same as for estimate (\ref{F})~in  Lemma \ref{JPh}.
 \end{proof}
 \noindent Now we shall  prove the strong convergence of the quantum random walks $p_t^{(h)}.$ Note that $J_t:\A \bigotimes \E(\K)\raro \A \bigotimes \Gamma$ 
 is the unique solution of the qsde
 \be \label{module-qsde}
 J_t =id_{\A\bigotimes \Gamma}+\int_0^t J_s  \Lambda_{\Theta}(ds).\ee 
  \noindent We define a family of maps $J_t^{(h)}$  by 
 \bean 
 && J_0^{(h)}(x \textbf{e}(f))u=x u \textbf{e}(f)\\
&& J_t^{(h)}(x \textbf{e}(f))u=J_{nh}^{(h)}(x \textbf{e}(f))u=x u \textbf{e}(f)+\sum_{k=1}^n
J_{(k-1)h} (\Lambda_{\Theta(x)}[k]\textbf{e}(f)) u \eean
for $t\in((n-1)h,nh].$
Thus by definition 
\be \label{j-tilda} J_t^{(h)}(x \textbf{e}(f))u=J_{nh}^{(h)}(x \textbf{e}(f))u=x u \textbf{e}(f)+\sum_{k=1}^n
 \Lambda_{j_{(k-1)h} (\Theta(x))}[k]u\textbf{e}(f).\ee
For $u\in {\mathbf h}_0, f\in \M$ the adapted process  $J_t$  satisfies 
\[J_t(x\textbf{e}(f))u =xu\textbf{e}(f)+\int_0^t J_s  \Lambda_{\Theta}(ds)  (x\textbf{e}(f)) u\] and the map  $t\mapsto J_t(x\textbf{e}(f))u$  is continuous. 
Thus by definition of this integral
   \[\lim_{h\raro 0}\|J_t(x\textbf{e}(f))u-J_t^{(h)}(x \textbf{e}(f))u\| =0\]
and hence
  \be \label{j-jhtilda} \lim_{h\raro 0}\|j_t(x)u \textbf{e}(f)-j_t^{(h)}(x)u \textbf{e}(f)\| =0.\ee
Now we are in position to prove the following  result
\bthm \label{st-con} Let $p_{t}^{(h)}$
be the quantum random walk  associated with $\beta(h).$
 Then for each $x\in \A $ and $ t\ge 0,  p_{t}^{(h)}(x)$  converges strongly to 
 $j_t(x).$ Thus $j_t:\A \raro \A\bigotimes \B(\Gamma)$ is a $*$-homomorphic flow.
 \ethm 
\begin{proof}
In order to  prove  
\be \label{ph-jh}\lim_{h\raro 0} \| p_t^{(h)}(x)u\textbf{e}(f)) - j_t(x)u\textbf{e}(f)\|=0,~ \forall  u\in {\mathbf{h}}_0, f\in \M,\ee
by (\ref{j-jhtilda}) it is sufficient to show that 
\be \label{jh-jhtilda}\lim_{h\raro 0} \| p_t^{(h)}(x)u\textbf{e}(f)) - j_t^{(h)}(x)u\textbf{e}(f)\|=0 \forall  u\in {\mathbf{h}}_0, f\in \M.\ee
For any fixed $ h>0, f\in \M$
let us define a family of bounded linear  maps
\[W_t^{(h)}:\A\rightarrow \A  {\mbox{$\bigotimes$}}  \Gamma\] given by, for $x\in \A$ and  $ u\in {\mathbf{h}}_0,$
\[W_t^{(h)}(x)u=p_t^{(h)}(x)u\textbf{e}(f)) - j_t^{(h)}(x)u\textbf{e}(f)\]
\[ =[{\mathcal P}_t^{(h)}(x\textbf{e}(f)) - J_t^{(h)}(x\textbf{e}(f)]u =:Y_t^{(h)}(x\textbf{e}(f))u. \]

\noindent Here, recall that  $\{J_t\}$ extend as a regular adapted process  mapping  $\A \bigotimes \B(\Gamma_{\rm fr}(\hat{\mathbf{k}}_0))\bigotimes \E(\C)$ into
$\A \bigotimes \B(\Gamma_{\rm fr}(\hat{\mathbf{k}}_0))\bigotimes \Gamma$ and  hence for each $ X\in \A \bigotimes \B(\Gamma_{\rm fr}(\hat{\mathbf{k}}_0)) $ the family $\{j_t(X)\}$ define by $j_t(X) \xi \textbf{e}(f)=J_t(X\otimes \textbf {e}(f))\xi, \forall \xi\in {\mathbf h}_0 \bigotimes \Gamma_{\rm fr}(\hat{\mathbf{k}}_0),f\in \C , $ is  a regular  $({\mathbf h}_0 \bigotimes \Gamma_{\rm fr},  \K)$-adapted process.
 For a given $f\in \M \subseteq \C$ by estimate (\ref{bigJtesti}), $ W_t^{(h)}$ extend as a bounded linear map from 
 $\A \bigotimes \B(\Gamma_{\rm fr}(\hat{\mathbf{k}}_0))$ into
$\A \bigotimes \B(\Gamma_{\rm fr}(\hat{\mathbf{k}}_0))\bigotimes \Gamma.$ 
 
Viewing $\A \bigotimes \B(\Gamma_{\rm fr}(\hat{\mathbf{k}}_0))$ and 
$\A \bigotimes \B(\Gamma_{\rm fr}(\hat{\mathbf{k}}_0))\bigotimes \Gamma$ as subspaces of $\B,$ let us denote by same symbol $W_t^{(h)}$ to  the canonical extentions of $W_t^{(h)}$  as linear maps from $\B$ into itself preserving the norm.

 In order to prove (\ref{jh-jhtilda}) we shall show that $\|W_t^{(h)}\|$  (as maps from $\B$ into itself) converges to $0$ as $h$ tends to $0.$ 
For any $X\in\A\bigotimes \B(\hat{\mathbf{k}}_0^{\mbox{\textcircled{\small m}}})$ and $ \xi \in {\mathbf h}_0\bigotimes
\hat{\mathbf{k}}_0^{\mbox{\textcircled{\small m}}}$
by (\ref{j-tilda}) and (\ref{qrw-mod}), we have
\bean \lefteqn{W_t^{(h)}(X)\xi }\\
&=&\sum_{k=1}^n\left [
N_{p_{(k-1)h}^{(h)}(\beta(h,X)-b(X))}[k] 
-\Lambda_{j_{(k-1)h} (\Theta(X))}[k]\right ]\xi \textbf{e}(f)\\
&&~~~-\sum_{k=1}^n p_{(k-1)h}^{(h)}  (X)(1-P_h[k])\xi \textbf{e}(f)\\  
&=&\sum_{k=1}^n \left[ \left (
N_{p_{(k-1)h}^{(h)}(\beta(h,X)-b(X))}[k]- N_{p_{(k-1)h}^{(h)}(\Theta(h,X))}[k]\right ) \right.
\xi \textbf{e}(f)\\
&&~~~~~~~~~+\left( N_{p_{(k-1)h}^{(h)}(\Theta(h,X))}[k] -\Lambda_{p_{(k-1)h}^{(h)} (\Theta(X))}[k]\right )\xi \textbf{e}(f)\\
&&~~~~~~~~~+\left.\left ( \Lambda_{p_{(k-1)h}^{(h)} (\Theta(X))}[k]
-\Lambda_{j_{(k-1)h} (\Theta(X))}[k]\right )\xi \textbf{e}(f)  \right]\\
&&-\sum_{k=1}^n p_{(k-1)h}^{(h)}  (X)(1-P_h[k])\xi \textbf{e}(f)
\eean
Using linearity of $N_{(\cdot)}[k]$ and $\Lambda_{(\cdot)}[k],$
\bean \lefteqn{\|W_t^{(h)}(X)\xi \|^2}\\
&&\le 4\left (  \| \sum_{k=1}^n 
   N_{p_{(k-1)h}^{(h)}(\beta(h,X)-b(X)
-\Theta(h,X))} [k]\xi_{k-1}\textbf{e}(f_{\left [ (k-1)h \right .})\|^2\right .\\
&&~~~+ \|\sum_{k=1}^n
 \left [ N_{p_{(k-1)h}^{(h)}(\Theta(h,X))}[k]-  \Lambda_{p_{(k-1)h}^{(h)}(\Theta(X))} [k] \right ]\xi_{k-1}\textbf{e}(f_{\left [ (k-1)h \right .})\|^2\\
 &&~~~+\|\sum_{k=1}^n p_{(k-1)h}^{(h)} (X)(1-P_h[k])\xi_{k-1}\textbf{e}(f_{\left [ (k-1)h \right .})\|^2 \\ 
&&~~~+\left .\|  \sum_{k=1}^n
  \Lambda_{\left [p_{(k-1)h}^{(h)} -j_{(k-1)h}\right ](\Theta(X))} [k] \xi_{k-1}\textbf{e}(f_{\left [ (k-1)h \right .})\|^2 \right )\\
&&=4(I_1+I_2+I_3+I_4).
\eean
 By Lemma \ref{NLamda}  we have
\[I_1+I_2+I_3\le  const(f,t) \|X\|^2 \|\xi\|^2 h.\]
 Now let us consider the terms in  $I_4.$ We have by estimate (3.13) in Proposition 3.3.5 in \cite{GS1}
\bean 
&&\|\sum_{k=1}^n
 \Lambda_{[p_{(k-1)h}^{(h)}-j_{(k-1)h}](\Theta(X))} [k] \xi\textbf{e}(f)\|^2\\
 &&=\|\sum_{k=1}^n Y_{(k-1)h}^{(h)} 
 \Lambda_{\Theta} [k] X\textbf{e}(f) \xi\|^2\\
 &&=\|\int_0^{nh} Y_s^{(h)} 
 \Lambda_{\Theta} (ds) X\textbf{e}(f) \xi\|^2\\
  &&\le  2e^t (1+\|f\|_\infty^2) \int_0^{nh} \|[(Y_s^{(h)}\otimes 1_{\hat{\mathbf{k}}_0})(\Theta(X)_{\hat{f}(s)}\textbf{e}(f))]
 \xi\|^2 ds.
 \eean
 It can be easily seen that 
 \[(Y_s^{(h)}\otimes 1_{\hat{\mathbf{k}}_0})(\Theta(X)_{\hat{f}(s)}\textbf{e}(f))
 =[(Y_s^{(h)}\otimes id_{\B(\hat{\mathbf{k}}_0)})(\Theta(X)\textbf{e}(f))]_{\hat{f}(s)}=[Y_s^{(h)}(\Theta(X)\textbf{e}(f))]_{\hat{f}(s)},\] so the above quantity is equal to
 \bean
 &&  2e^t (1+\|f\|_\infty^2) \int_0^{nh} \|[Y_s^{(h)}(\Theta(X)\textbf{e}(f))]_{\hat{f}(s)
 }\xi\|^2 ds,\\
 &&=  2e^t (1+\|f\|_\infty^2) \int_0^{nh} \|W_s^{(h)}(\Theta(X)) \xi\otimes \hat{f}(s)\|^2 ds\\
 &&\le 2e^t (1+\|f\|_\infty^2)^2  \sum_{k=1}^n h \|W_{(k-1)h}^{(h)}\|^2 \|\Theta(X))\|^2 \|\xi \|^2 \\
 &&\le c_f \sum_{k=1}^n h \|W_{(k-1)h}^{(h)} \|^2\|\Theta\|^2 \|X\|^2 \|\xi \|^2 .
 \eean
 \noindent Combining all the above estimates, we obtained 
 \bea \label{Wt}\lefteqn{\|W_t^{(h)}(X)\xi \|^2}\\
 &\le &h C \|X\|^2 \|\xi \|^2+ D \sum_{k=1}^n h
\|W_{(k-1)h}^{(h)} \|^2 \|X\|^2 \|\xi \|^2 \nonumber
\eea
for some constant $C$ and $D$ independent of $h.$ 
For any $X\in \B$ and $\xi\in{\mathbf h}_0 \bigotimes \Gamma_{\rm{fr}}(\hat{\mathbf{k}}_0) \bigotimes \Gamma$
we can write $ X=\bigoplus_{m\ge 0} X_m\oplus X^\prime$ with $X_m\in \A\bigotimes \B(\hat{\mathbf{k}}_0^{\mbox{\textcircled{\small m}}}),X^\prime \in \B_c $ and $\xi=\bigoplus_{m\ge 0} \xi_m\oplus \xi^\prime$ with 
$\xi_m\in {\mathbf h}_0 \bigotimes \hat{\mathbf{k}}_0^{\mbox{\textcircled{\small m}}}$ and $\xi^\prime$ belong to orthogonal complement of ${\mathbf h}_0 \bigotimes \hat{\mathbf{k}}_0^{\mbox{\textcircled{\small m}}}$ for all $m\ge 0.$
Using the estimate (\ref{Wt}) we have
\bean \lefteqn{\|W_t^{(h)}(X)\xi \|^2}\\
&=&\sum_{m\ge 0}\|W_t^{(h)}(X_m)\xi_m \|^2\\
 &\le &h C \sum_{m\ge 0}\|X_m\|^2 \|\xi_m \|^2+ D \sum_{k=1}^n h
\|W_{(k-1)h}^{(h)} \|^2 \sum_{m\ge 0}\|X_m\|^2 \|\xi_m \|^2\\ 
&\le &h C \|X\|^2 \|\xi \|^2+ D \sum_{k=1}^n h
\|W_{(k-1)h}^{(h)} \|^2 \|X\|^2 \|\xi \|^2   
\eean
 Taking supremum over all $\xi\in  {\mathbf h}_0 \bigotimes \Gamma_{\rm{fr}}(\hat{\mathbf{k}}_0) \bigotimes \Gamma, X\in \B$ such that $\|\xi \|\le 1,\|X\|\le  1$ we get
\be \label{w-cb}\|W_t^{(h)}\|^2=\|W_{nh}^{(h)}\|^2\le h C +h D \sum_{k=1}^n
\|W_{(k-1)h}^{(h)}\|^2  .\ee
By definition  $\|W_0^{(h)}\|^2=0$ so (\ref{w-cb}) gives
$\|W_{h}^{(h)}\|^2 \le hc$ and
\[\|W_{2h}^{(h)}\|^2 \le hc+hD\|W_{h}^{(h)}\|^2 \le ch(1+hD).\]
Then by induction it follows that 
 \[\|W_t^{(h)}\|^2=\|W_{nh}^{(h)}\|^2\le h C (1+  h D)^{n-1} \le h Ce^{Dt} \]
and hence  
\[\lim _{h\rightarrow 0}\|W_t^{(h)}\|^2=0, ~\mbox{
inparticular  }~\lim_{h\raro 0} \| p_t^{(h)}(x)u\textbf{e}(f)) - j_t^{(h)}(x)u\textbf{e}(f)\|=0 .\]
Which says that for any  $u\in {\mathbf{h}}_0$ and $f\in \M, ~\{p_t^{(h)}(x)u\textbf{e}(f)):h>0\}$ is  Cauchy in ${\mathbf h}_0 \bigotimes \Gamma.$ Since  $\|p_t^{(h)}(x)\|\le \|x\| $ and algebraic tensor product ${\mathbf{h}}_0 \bigotimes \E(\M)$ is dense in ${\mathbf{h}}_0 \bigotimes \Gamma$ it follows that  $\{p_t^{(h)}(x) \xi:h>0\}$ is Cauchy  for all  $\xi\in {\mathbf{h}}_0 \bigotimes \Gamma$ and hence for each $x\in \A, ~\{p_t^{(h)}(x)\}$
converges strongly to $j_t(x).$ Thus $j_t:\A \raro \A \bigotimes \B(\Gamma)$  is a contractive $*$-homomorphic flow.
\end{proof}

\brmrk (i) It may be observed that in the above quantum stochastic dilation  $\{j_t\}$ of the dynamical semigroup $\{T_t\}$ there is no ``Poisson" term since $\theta_4(x)=0$ for all $x\in \A.$ This is only to be expected since the choice of representation of $\A$  is $x\otimes 1_{{\mathbf k}_0}$ for all $x\in \A.$
The more general case of dilation using the convergence of quantum random walks where the representation is non trivial (and  therefore  will have non zero ``Poisson"  component) is being investigated.\\
(ii)  The method of proof employed above does not seem to be amenable  to adaptation for a dynamical semigroup  with unbounded generator. On the other hand,  one has example of  the convergence of random walks to diffusion processes
(which of course, has unbounded generators )  in the  classical case. For the handling of these cases, one may have to find  different method to replace the proof of Theorem  \ref{st-con}.
\ermrk

\end{document}